\documentclass[11pt]{article}
\usepackage{amsfonts}
\usepackage{dsfont}
\usepackage{amsfonts}
\usepackage{bbm}
\usepackage{mathrsfs}
\usepackage{color}
\usepackage{hyperref}
\usepackage{mathrsfs}
\usepackage{amsmath}
\usepackage{enumitem}
\allowdisplaybreaks[4]
\usepackage{amsfonts}
\usepackage{amsthm}
\usepackage{amssymb, amsmath, cite}
\usepackage{bm}
\usepackage{amsthm}
\usepackage{verbatim}
\usepackage{graphicx}
\usepackage{color}
\usepackage{amsfonts}
\usepackage{dsfont}
\usepackage{amsfonts}
\usepackage{bbm}
\usepackage{mathrsfs}
\usepackage{color}
\usepackage{graphicx,amsmath,amsfonts,amssymb,amscd,amsthm}
\usepackage{hyperref}
\usepackage[numbers,sort&compress]{natbib}
\usepackage{enumitem}
\setenumerate[1]{itemsep=1pt,partopsep=0pt,parsep=\parskip,topsep=1pt}
\setitemize[1]{itemsep=1pt,partopsep=0pt,parsep=\parskip,topsep=1pt}

\newtheorem{Theorem}{Theorem}[section]
\newtheorem{Lemma}{Lemma}[section]

\newtheorem{Remark}{Remark}[section]
\newtheorem{Definition}{Definition}[section]
\newtheorem{Corollary}{Corollary}

\usepackage{graphicx,amsmath,amsfonts,amssymb,amscd,amsthm}
\usepackage[numbers,sort&compress]{natbib}
\usepackage{array}
\usepackage{booktabs}
\usepackage{multirow}
\usepackage{makecell}
\usepackage{abstract}
\hypersetup{colorlinks,linkcolor=blue,citecolor=blue}
\newtheoremstyle{kai}
{3pt} {3pt} {} {} {\bfseries} {.} {.5em} {}
\def\EquationsBySection{\def\theequation
{\thesection.\arabic{equation}}
\@addtoreset{equation}{section}}

\newcommand\old[1]{}

\newcommand{\Bp}{\begin{proof}}
 \newcommand{\Ep}{\end{proof}}
\renewcommand{\theequation}{\thesection.\arabic{equation}}

\numberwithin{equation}{section}
\addcontentsline{toc}{section}{Acknowledgement}

\newcommand{\newcom}{\newcommand}
\newcom{\R}{\mathbb R}
\newcom{\N}{\mathbb N}
\newcom{\e}{\varepsilon}

\newcom{\al}{\alpha}
\newcom{\be}{\beta}
\newcom{\del}{\delta}

\newcom{\ga}{\gamma}
\newcom{\Ga}{\Gamma}

\newcom{\Lam}{\Lambda}
\newcom{\lam}{\lambda}

\newcom{\Om}{\Omega}
\newcom{\om}{\omega}

\newcom{\Si}{\Sigma}
\newcom{\si}{\sigma}
\newcom{\s}{\varsigma}

\newcom{\tht}{\theta}
\newcom{\dtri}{\nabla}
\newcom{\tri}{\triangle}
\newcom{\mf}{\,\,\,\,\mathrm{for\,\,all} \,\,}
\newcom{\qmf}{\qquad\mathrm{for\,\,all}~~ }

\newcom{\no}{\nonumber}
\newcom{\f}{\frac}
\newcom{\na}{\nabla}
\newcom{\Del}{\Delta}
\newcom{\ep}{\epsilon}

\newcom{\p}{\partial}

\newcom{\uep}{u_{\e}}
\newcom{\vep}{v_{\e}}
\newcom{\nep}{n_{\epsilon}}
\newcom{\cep}{{c}_{\epsilon}}
\newcom{\whc}{\widehat{c}}
\newcom{\ocep}{\widehat{c}_{\epsilon}}
\newcom{\cp}{{c}_{\epsilon}}
\newcom{\iom}{\int_{\Omega}}
\newcom{\wcep}{\widetilde{c}_{\e}}
\newcom{\ce}{c^*}

\newcom{\beq}{\begin{equation}}
\newcom{\eeq}{\end{equation}}
\newcom{\beno}{\begin{eqnarray*}}
\newcom{\eeno}{\end{eqnarray*}}

\begin{document}
\title{\bf Refined existence theorems for doubly degenerate chemotaxis-consumption systems with large initial data}

\author{{\sc Duan Wu\thanks{E-mail: duan\underline{~}wu@126.com}}\\
{\small Institut f\"ur Mathematik, Universit\"at Paderborn, 33098 Paderborn, Germany}\\
{\small and}\\
{\small School of Mathematics and Statistics, Northwestern Polytechnical University, Xi'an 710129, PR China}}
\date{}
\maketitle
\begin{abstract}
 This work considers the doubly degenerate nutrient model
\begin{equation*}\label{AH1}
   \left\{
   \begin{split}
     &u_t=\nabla\cdot\left(u^{m-1}v\nabla u\right)-\nabla\cdot\left(f(u)v\nabla v\right)+\ell uv,&&x\in\Omega,\,t>0,\\
     &v_t=\Delta v-uv, &&x\in\Omega,\,t>0,
   \end{split}
   \right.
\end{equation*}
under no-flux boundary conditions in a smoothly bounded convex domain $\Omega\subset \R^n$ ($n\le 2$), where the nonnegative function $f\in C^1([0,\infty))$ is assumed to satisfy $f(s)\le C_fs^{\al}$ with $\al>0$ and $C_f>0$ for all $s\ge 1$.

When $m=2$, it was shown that a global weak solution exists, either in one-dimensional setting with $\al=2$, or in two-dimensional version with $\al\in(1,\f{3}{2})$. The main results in this paper assert the global existence of weak solutions for $1\le m<3$ and classical solutions for $3\le m<4$ to the above system under the assumption
\begin{equation*}
  \al\in \left\{
   \begin{split}
     &\left[m-1,\min\left\{m,\f{m}{2}+1\right\}\right]~~&&\textrm{if}~~n=1,\quad\quad\textrm{and}\\
     &\left(m-1,\min\left\{m,\f{m}{2}+1\right\}\right)~~&&\textrm{if}~~n=2,
   \end{split}
   \right.
\end{equation*}
which extend the range $\al\in(1,\f{3}{2})$ to $\al\in(1,2)$ in two dimensions for the case $m=2$. Our proof will be based on a new observation on the coupled energy-type functional and on an inequality with general form.
\textbf{AMS (2020) Subject Classification}: {35K65; 92C17; 35A01}
\end{abstract}

\section{Introduction}
Being motivated by some experimental indications when investigating the dynamical motility of \textit{Bacillus subtilis} grown on thin agar plates (\cite{fujikawa1992physica},\cite{m-fujikawa1990physica},\cite{ohgiwari-m-m1992physica}), Kawasaki et.al in \cite{kawasaki-mmus1997jtb} introduced the following reaction-diffusion continuous system
\begin{equation}\label{0s}
\left\{
\begin{split}
&u_t=\nabla\cdot\left(D(u,v)\nabla u\right)-\nabla\cdot\left(S(u,v)\nabla v\right)+\ell uv,\\
&v_t=\Delta v-uv
\end{split}
\right.
\end{equation}
with the non-linear diffusion coefficient $D(u,v)=uv$ in line with the peculiarity that the bacteria would be immobile when exposed to regions with the low level of nutrients, and with the absence of taxis processes given by $S(u,v)=0$. Subsequently, to comply with experimental phenomena discovered in \cite{golding.etal1998} and \cite{ben-jacob.etal2000}, Leyva et.al in \cite{Leyva.etal2013physica} proposed that the bacterial response function in the chemotactic flux term should be considered by $S(u,v)=\chi u^2v$ with $\chi>0$ measuring the intensity of the chemotaxis, which was supported by numerical simulations: the conformity between the experimentally gained observations and the behavior of numerical solutions is noticeably greater than the case taking $S(u,v)=0$ (see \cite[Sections 4 and 5]{Leyva.etal2013physica}).

Within the realm of results from mathematical analysis, even in the case of the prototypical setting upon the choices $D(u,v)=1$ and $S(u,v)=u$, to the best of our knowledge there has been sparse literature with respect to solvability and spatially homogeneous equilibria for \eqref{0s} (\cite{tao-winkler2012jde}, \cite{winkler2017jde}), which are confined to the relatively lower-dimensional space.

For the doubly degenerate diffusion problem such as when $D(u,v)=uv$, virtually no mathematical findings regarding basic solution theory touched on it until Winkler in \cite{winkler2021tams} revealed the existence of global weak solutions for arbitrary large initial data in one-dimensional space with the choice of $S(u,v)=u^2v$ as suggested in \cite{Leyva.etal2013physica}, where the most stunning feature is that the first component of the solution asymptotically stabilizes towards a nontrivial function coinciding with the spatial profile of a solution to a scalar porous medium-type parabolic equation in stark contrast to the majority of previous results concerning the large time behavior. In the corresponding taxis-free framework, namely $S(u,v)=0$, the global existence of weak solutions and the stabilization were asserted in convex domains of any dimension (\cite{winkler2022cvpde}). For the same choices of $S(u,v)$ and $D(u,v)$ as studied in \cite{winkler2021tams}, Li and Winkler removed the integrability assumption $\iom \ln u_0>-\infty$ to establish the similar results to system \eqref{0s} in the spatially one-dimensional analogue by relying on a striking energy functional (\cite{li-winkler2022cpaa}), while for its two-dimensional version the global solvability were substantiated in bounded convex domains with a smallness condition exclusively involving $v_0$ (\cite{winkler2022na}). When the cross-diffusion coefficient has the form $S(u,v)=\chi u^{\al}v$ with $\al>0$, without any smallness restrictions on the size of initial data in bounded convex domains, Li identified that \eqref{0s} possesses a global weak solution under the range $1<\al<\f{3}{2}$ in two-dimensional case, and $\f{7}{6}<\al<\f{13}{9}$ in three-dimensional case (\cite{li2022jde}).

Motivated by the above precedent results, we have a natural question that when $D(u,v)=uv$ and $S(u,v)=\chi u^{\al}v$ in \eqref{0s}, if $\al=\f{3}{2}$ is a critical existence exponent in two dimensions. Specifically, we attempt to figure out whether the smallness on initial data is necessary to guarantee the solvability in the flavor of \cite{li2022jde} for the case $\f{3}{2}\le \al<2$. Meanwhile, we are devoted to studying the corresponding solvability for system \eqref{0s} with a more general nonlinear diffusion term, which plays an important role in the biological context (see \cite{kit1997}, \cite{Bellomo2022MMM}).

More precisely, in the present work we will be concerned with the initial-boundary value problem
\begin{equation}\label{s}
\left\{
\begin{split}
&u_t=\nabla\cdot\left(u^{m-1}v\nabla u\right)-\nabla\cdot\left(f(u)v\nabla v\right)+\ell uv,&&x\in\Omega,\,t>0,\\
     &v_t=\Delta v-uv, &&x\in\Omega,\,t>0,\\
&\left(u^{m-1}v\nabla u-f(u)v\nabla v\right)\cdot\nu=\na v\cdot\nu=0,&& x\in\partial\Omega,\,t>0,\\
&u(x,0)=u_0(x),~~v(x,0)=v_0(x),&&x\in\Omega,
\end{split}
\right.
\end{equation}
with $m\ge 1$, $\ell\ge 0$, under the overall hypotheses that $f$ besides
\beq\label{f}
f\in C^1([0,\infty))~~{\rm{and}}~~ f(u)\ge 0\qmf u\ge0
\eeq
is such that
\beq\label{f1}
f(u)\le C_fu(u+1)^{\al-1} \qmf u\ge0
\eeq
or
\beq\label{f2}
f(u)\le C_fu^{\al} \qmf u\ge0
\eeq
with $C_f>0$ and $\al\ge 0$. In addition, the initial data $(u_0, v_0)$ is throughout supposed to satisfy that
\begin{equation}\label{indata}
\left\{
\begin{split}
&u_0\in W^{1,\infty}(\Om){\rm{~is~nonnegative~with~}} u_0\not\equiv 0\qquad{\rm{and}}\\
&v_0\in W^{1,\infty}(\Om) {\rm{~is~positive~in~}}\overline{\Om}.
\end{split}
\right.
\end{equation}
Before stating our main results, we first introduce the definition of weak solutions pursued in this context.
\begin{Definition}\label{de}
Let $n\in\{1,2\}$ and $\Omega\subset \mathbb{R}^n$ be a bounded domain with smooth boundary, $m\ge 1$ and $\ell\ge 0$. Suppose that $f$ satisfies \eqref{f1} or \eqref{f2} with $C_f>0$ and $\al\ge0$ in addition to \eqref{f}, and that $u_0\in L^1(\Om)$ and $v_0\in L^1(\Om)$ are nonnegative. Then we call that a pair $(u,v)$ of nonnegative functions
\beno\label{de1}
\left\{
\begin{split}
&u\in L_{loc}^{1}(\overline\Om\times[0,\infty))\qquad{\rm{and}}\\
&v\in L_{loc}^{\infty}(\overline\Om\times[0,\infty))\cap L_{loc}^{1}([0,\infty);W^{1,1}(\Om)).
\end{split}
\right.
\eeno
satisfying
\beno\label{de2}
u^m\na v\in L_{loc}^{1}\left(\overline\Om\times[0,\infty);\R^n\right) \quad{\rm{and}}\quad
u^mv\in L_{loc}^{1}\left(\overline\Om\times[0,\infty)\right)
\eeno
is a global weak solution of \eqref{s} if
\begin{align*}
-\int_{0}^{\infty}\iom u\varphi_t-\iom u_0\varphi(\cdot,0)&=\f{1}{m}\int_{0}^{\infty}\iom u^m\na v\cdot\na\varphi+\f{1}{m}\int_{0}^{\infty}\iom u^mv\Delta \varphi\no\\
&\quad+\int_{0}^{\infty}\iom f(u)v\na v\cdot\na\varphi+\ell\int_{0}^{\infty}\iom uv\varphi
\end{align*}
for all $\varphi\in C_{0}^{\infty}\left(\overline\Om\times[0,\infty)\right)$ fulfilling $\f{\partial\varphi}{\partial\nu}=0$ on $\partial\Om\times(0,\infty)$, as well as
\beno\label{de4}
\int_{0}^{\infty}\iom v\varphi_t+\iom v_0\varphi(\cdot,0)=\int_{0}^{\infty}\iom \na v\cdot\na\varphi+\int_{0}^{\infty}\iom uv\varphi
\eeno
for each $\varphi\in C_{0}^{\infty}\left(\overline\Om\times[0,\infty)\right)$.
\end{Definition}

The first result concerns one-dimensional space.
\begin{Theorem}\label{th1}
Let $\Omega\subset \R$ be an open bounded interval. Suppose that $\ell\ge0$ and that the initial data $u_0$ and $v_0$ satisfy \eqref{indata}. Then for all $K>0$ with the property that
\beq\label{indata1}
\|u_0\|_{L^{\infty}(\Om)}+\|v_0\|_{L^{\infty}(\Om)}+\|\partial_x \ln v_0\|_{L^{\infty}(\Om)}\le K,
\eeq
if one of the following cases holds:\\
(i) $1\le m<2$, $f$ fulfills \eqref{f} and \eqref{f1} with $m-1\le\al\le m$;\\
(ii) $2\le m<3$, $f$ fulfills \eqref{f} and \eqref{f2} with $m-1\le\al\le \f{m}{2}+1$,\\
the problem \eqref{s} admits a global weak solution $(u,v)$ in the sense of Definition \ref{de}. Moreover, if\\
(iii) $3\le m<4$, $f$ fulfills \eqref{f} and \eqref{f2} with $m-1\le\al\le \f{m}{2}+1$ and $u_0>0$ in $\overline{\Om}$,\\
the problem \eqref{s} admits a global classical solution $(u,v)$.

Furthermore, $v>0$ a.e. in $\Om\times(0,\infty)$, and that there exists $C(K)>0$ such that
\beq\label{th1r}
{\|u(\cdot,t)\|_{L^{\infty}(\Om)}}+\|v(\cdot,t)\|_{W^{1,\infty}(\Om)}\le C(K) {{\qquad for~ a.e.\quad}} t>0.
\eeq
\end{Theorem}
\begin{Remark}
 {\rm{In addition, with the aid of the Harnack type inequality valid in one-dimensional space (\cite{winkler-preprint}), in much the same way as used in \cite{li-winkler2022cpaa}, we could further obtain the following asymptotic behavior result.}}
\begin{Corollary}
Let the assumptions in Theorem \ref{th1} hold and $(u,v)$ be as accordingly given by Theorem \ref{th1}, then there exists $u_{\infty}=w(\cdot,1)\in C^0(\overline\Om)$ such that $u(\cdot,t)\rightarrow u_{\infty}$ in $L^{\infty}(\Om)$ as $t\rightarrow\infty$. Here $w\in C^0{(\overline\Om\times[0,1])}\cap L_{loc}^2([0,1];W^{1,2}(\Om))$ is a weak solution of
\beno
\left\{
\begin{split}
&w_{\tau}=\big(a(x,\tau)w^{m-1}w_x\big)_x-\big(b(x,\tau)f(w)\big)_x+\ell a(x,\tau)w,&&x\in\Omega,\,\tau\in(0,1),\\
&w_{x}=0, &&x\in\p\Omega,\,\tau\in(0,1),\\
&w(x,0)=u_0(x), &&x\in\Omega,\\
\end{split}
\right.
\eeno
where $a(x,\tau)$ and $b(x,\tau)$ are defined as that in \cite[Theorem 1.3]{winkler2021tams}.
\end{Corollary}
\end{Remark}

In two-dimensional setting, for convex domains, we have the following result.
\begin{Theorem}\label{th2}
Let $\Omega\subset \R^2$ be a bounded convex domain with smooth boundary, and let $\ell\ge0$. Suppose that the initial data $u_0$ and $v_0$ satisfy \eqref{indata}. Then for all $K>0$ with the property that
\beq\label{indata2}
\|u_0\|_{L^{\infty}(\Om)}+\|v_0\|_{L^{\infty}(\Om)}+\|\na \ln v_0\|_{L^{\infty}(\Om)}\le K,
\eeq
if one of the following cases holds:\\
(i) $1\le m<2$, $f$ fulfills \eqref{f} and \eqref{f1} with $m-1<\al< m$;\\
(ii) $2\le m<3$, $f$ fulfills \eqref{f} and \eqref{f2} with $m-1<\al< \f{m}{2}+1$,\\
the problem \eqref{s} admits a global weak solution $(u,v)$ in the sense of Definition \ref{de}. Moreover, if\\
(iii) $3\le m<4$, $f$ fulfills \eqref{f} and \eqref{f2} with $m-1<\al<\f{m}{2}+1$ and $u_0>0$ in $\overline{\Om}$,\\
the problem \eqref{s} admits a global classical solution $(u,v)$.

Furthermore, $v>0$ a.e. in $\Om\times(0,\infty)$, and that for all $p>2$, there exists $C(p,K)>0$ such that
\beq\label{th2r}
\|u(\cdot,t)\|_{L^p(\Om)}+\|v(\cdot,t)\|_{W^{1,\infty}(\Om)}\le C(p,K) {{\qquad for~ a.e.\quad}} t>0.
\eeq
\end{Theorem}

\begin{Remark}
{\rm{The large time behavior has not yet been solved because of the absence of the corresponding Harnack type inequality in two-dimensional space, which also results in that we can only obtain the $L^p$ bounds for $u$ in \eqref{th2r} instead of the $L^{\infty}$ bounds.}}
\end{Remark}

\begin{Remark}
{\rm{When $m=2$, the upper {bound} on $\al$ obtained in this paper seems natural and reasonable from the perspective that if $\al=2$, the solvability requires imposing small conditions on the initial data as indicated in \cite{winkler2022na}. }}
\end{Remark}

\begin{Remark}
{\rm{In a manner rather similar to proving Theorem \ref{th2}, the corresponding results in higher-dimensional space is able to be derived: at least in the context of three dimensions with the case $m=2$, we can check that the global solvability holds within the range $1<\al<\f{3}{2}$, which together with Theorem \ref{th2} extends the previous results in \cite{li2022jde}.}}
\end{Remark}
The foundation of our approach relies closely on the useful functional first discovered in \cite{winkler2022dcdsb} of the form
\beno
\iom \f{|\na v|^{q}}{v^{q-1}},~~q\ge 2
\eeno
with the favorable features
\begin{align*}
&\f{d}{dt}\iom \f{|\na v|^{q}}{v^{q-1}}+\f{q}{2}\iom \f{|\na v|^{q-2}}{v^{q-3}}|D^2\ln v|^2+(q-1)^2\iom u\f{|\na v|^q}{v^{q-1}}\le C_1(q)\iom u^{\f{q+2}{2}}v,
\end{align*}
in particular when $q=4$, by a slightly different estimation used to get the above estimation and making use of the good term $\iom \f{|\na v|^{6}}{v^{5}}$ (identified by Lemma \ref{lephi}),
\begin{align}\label{02}
&\f{d}{dt}\iom \f{|\na v|^{4}}{v^{3}}+\iom \f{|\na v|^{2}}{v}|D^2\ln v|^2+\iom u\f{|\na v|^4}{v^3}\le C_2\iom v|\na u|^2,
\end{align}
where $C_1(q)$ and $C_2$ are positive constants (Lemma \ref{lenavq}).

In the spatially one-dimensional case, we commence with a coupled functional expressed by
\beq\label{0F}
\mathcal{F}(t):=\iom u^p+\iom \f{v_x^{2}}{v},
\eeq
where $p$ is suitably taken by $p=3-m$, roughly speaking. For the peculiar cases of $m=2$ and $m=3$, we substitute $1$ with the function $\ln u$ which {exhibits} the same scaling as $1$, and correspondingly $u$ is replaced by $u\ln u$ (Lemma \ref{leF}). Then on the basis of $L^4$ bound on $v_x$ (Lemma \ref{lenav4}), by means of a functional inequality (Lemma \ref{lefi1}), the ill-signed contributions of $\mathcal{F}$ could be controlled by the dissipated integral $\iom u\f{v_x^2}{v}$ and the favorable term $\iom vu_x^2$ produced by the evolution of first summand in \eqref{0F}. This furnishes some boundedness properties presented in Lemma \ref{lecru}, which {will} be utilized to derive bounds for $u$ in $L^p$ with any $p>1$ (Lemma \ref{leup5}).

When the spatial setting is two-dimensional, it is noticeable from Lemma \ref{leup1} that if we have a bound for $u$ in $L^{p_0}$ with some $p_0>1$, then {using an} energy-like functional of the form
\beq\label{0G0}
\iom u^p+\iom \f{|\na v|^{q}}{v^{q-1}}
\eeq
with suitably chosen $q\in \left(\f{2(p+m-1)}{p_0},~2(p_0+p+m-2)\right)=:I$ is an effectively direct method to obtain the $L^{p}$ bounds for $u$ with arbitrary $p>2$. This way is not applicable when we try to merely rely on the basic boundedness of $\|u\|_{L^1}$ because $I$ is empty when $p_0=1$. Thus, the current difficulty consists in how to obtain an $L^{p_0}$ bound for $u$ with some $p_0>1$. To achieve this, the core of our analysis is based on the following quasi-energy functional, given by
\beq\label{0G}
\mathcal{G}(t):=c\iom u^p+\iom \f{|\na v|^{4}}{v^3}+\iom u^{p^*},
\eeq
where $c>0$, $p$ is as in \eqref{0F} and $p^*\in (2-m,3-m)$. Some basic differential inequalities for the first two summands in \eqref{0G} are shown in Lemma \ref{leG}, while the key role in addressing the unfavorable expressions therein is played by the crucial observation that
\beno\label{0u-v}
\iom u^{\beta}v\le \eta\iom u^{\kappa}v|\na u|^2+\eta\iom \f{|\na v|^{2}}{v}|D^2\ln v|^2+C_3(\beta,\eta)\iom uv
\eeno
holds for $\beta\in[1,\kappa+3)$ with $\kappa\in(-1,0)$ (Lemmas \ref{leu-v}). The integral $\iom \f{|\na v|^{2}}{v}|D^2\ln v|^2$ can be dissipated during the evolution of $\iom \f{|\na v|^{4}}{v^{3}}$ as shown in \eqref{02}, and the other integral $\iom u^{\kappa}v|\na u|^2$ will be absorbed by further testing $u^{p^*}$ in a standard procedure with diffusion-induced quantity $\iom u^{p^*+m-3}v|\na u|^2$ arising. Then we could get some essential boundedness properties, which serve as a crucial ingredient in our derivation of an $L^{p_0}$ bound for $u$ (Lemmas \ref{lecru2-1}-\ref{lecru2-3}).

Based on these results, we are able to collect some higher regularities (cf. Lemma \ref{lenav-inf}, Lemma \ref{le-u-infty} and  Lemma \ref{lest-1}) and further announce the conclusion on the global existence.

\section{Some preliminaries}
In this section, we introduce some basic lemmas. Similar with the approximating procedure used in \cite{li-winkler2022cpaa}, we consider the regularized variant of \eqref{s} given by
\begin{equation}\label{s1}
\left\{
\begin{split}
&{\uep}_t=\nabla\cdot\left(u_{\e}^{m-1}{\vep}\nabla {\uep}\right)-\nabla\cdot\left(f({\uep}){\vep}\nabla {\vep}\right)+\ell {\uep}{\vep},&&x\in\Omega,\,t>0,\\
&{\vep}_t=\Delta {\vep}-{\uep}{\vep}, &&x\in\Omega,\,t>0,\\
&\f{\p\uep}{\p\nu}=\f{\p\vep}{\p\nu}=0,&& x\in\partial\Omega,\,t>0,\\
&\uep(x,0)=u_{0\e}(x),~~v(x,0)=v_{0\e}(x):=v_0(x),&&x\in\Omega
\end{split}
\right.
\end{equation}
with $\e\in(0,1)$, where {$u_{0\e}(x)$ depending on $m$ is defined by
\begin{equation}\label{u0e}
u_{0\e}(x):=\left\{
\begin{split}
&u_0(x)+\e, &&1\le m<3,\\
&u_0(x), &&3\le m<4.\\
\end{split}
\right.
\end{equation}}
In light of the positivity of $u_{0\e}$, this problem enjoys the following properties.

\begin{Lemma}\label{lelocal}
Suppose that the assumptions in Theorem \ref{th1} {or} Theorem \ref{th2} are satisfied. Then for each $\e\in(0,1)$, there exist $T_{\max,\e}\in(0,\infty]$ and at least one pair $(\uep,\vep)$ of functions
\beq\label{local}
\left\{
\begin{split}
&\uep\in \cap_{q\ge1} C^0\left([0,T_{\max,\e});W^{1,q}(\Omega)\right)\cap C^{2,1}\left(\overline\Omega\times (0,T_{\max,\e})\right)\\
&\vep\in \cap_{q\ge1} C^0\left([0,T_{\max,\e});W^{1,q}(\Omega)\right)\cap C^{2,1}\left(\overline\Omega\times (0,T_{\max,\e})\right)
\end{split}
\right.
\eeq
such that $\uep>0$ and $\vep>0$ in $\overline\Omega\times (0,T_{\max,\e})$, that $(\uep,\vep)$ solves \eqref{s1} in the classical sense in $\Om\times (0,T_{\max,\e})$, and that
\beq\label{exten}
 if~T_{max,\e}<\infty,~~then~~
\limsup_{t\rightarrow T_{max,\e}}\|\uep(\cdot, t)\|_{L^{\infty}(\Omega)}=\infty.
\eeq
In addition, this solution satisfies
\beq\label{vin}
\|\vep(\cdot,t)\|_{L^{\infty}({\Om})}\leq \|\vep(\cdot,t_0)\|_{L^{\infty}({\Om})}~~~~for~all~t_0\in[0,T_{max,\e})~and~any~t\in(t_0,T_{max,\e}),
\eeq
\beq\label{u1}
\iom u_{0\e}\le\iom \uep(\cdot,t)\le\iom u_{0\e}+\ell\iom v_{0\e}~~~~for~all~t\in(0,T_{max,\e})
\eeq
and
\beq\label{uv1}
\int_{t_0}^{T_{max,\e}}\iom \uep\vep\leq \iom \vep(\cdot,t_0)~~~~for~all~t_0\in[0,T_{max,\e}).
\eeq
\end{Lemma}
\Bp
By means of Theorem 14.4, Theorem 14.6 and Theorem 15.5 in \cite{Amann}, we can obtain the local existence and extensibility with the form that
\begin{multline}\label{exten-1}
{\rm{if}}~T_{max,\e}<\infty,~~{\rm{then}}\\
\limsup_{t\rightarrow T_{max,\e}}\bigg\{\|\uep(\cdot, t)\|_{W^{1,\infty}(\Omega)}+\|\vep(\cdot, t)\|_{W^{1,\infty}(\Omega)}+\bigg\|\f{1}{\uep(\cdot, t)}\bigg\|_{L^{\infty}(\Omega)}+\bigg\|\f{1}{\vep(\cdot, t)}\bigg\|_{L^{\infty}(\Omega)}\bigg\}=\infty.
\end{multline}
We assert that \eqref{exten-1} and \eqref{exten} are equivalent within this framework. It is obvious that \eqref{exten} implies \eqref{exten-1}. On the other hand, after assuming that for some $\e\in(0,1)$ there exists $c_1(\e)>0$ such that $\|\uep(\cdot, t)\|_{L^{\infty}(\Omega)}\le c_1(\e)$ for all $t\in(0,T_{max,\e})$ with $T_{max,\e}<\infty$, we could find positive constants $c_2(\e)$ and $c_3(\e)$ such that
\beq\label{0vnav-inf}
\|\vep(\cdot, t)\|_{W^{1,\infty}(\Omega)}\le c_2(\e)\qmf t\in(0,T_{max,\e})
\eeq
and
\beq\label{0v-up}
\vep\ge c_3(\e) \qquad {\rm{in}}~~\Om\times(0,T_{max,\e})
\eeq
by using the same way as the first part of \cite[Lemma 2.1]{winkler2022na}. It is easy to verify that there exists $c_4(\e)>0$ such that
\beq\label{fu}
f(\uep)\le c_4(\e)u_{\e}^{\f{m-1}{2}}
\eeq
either under the case $1\le m<2$, $f$ fulfills \eqref{f} and \eqref{f1} with $m-1\le\al\le m$ or the case $2\le m<4$, $f$ fulfills \eqref{f} and \eqref{f2} with $m-1\le\al\le\f{m}{2}+1$. Thus, rewriting the first equation in \eqref{s1} in the following form
\beno
{\uep}_t=\nabla\cdot a_{\e}(x,t,\uep,\na\uep)+b_{\e}(x,t,\uep,),\qquad x\in\Omega,\,t\in(0,T_{max,\e})
\eeno
with
\begin{align*}
&a_{\e}(x,t,\uep,\na\uep):=\vep(x,t)\uep^{m-1}(x,t)\na\uep(x,t)-f(\uep(x,t))\vep(x,t)
\na\vep(x,t)\qquad {\rm{and}}\no\\
&b_{\e}(x,t)=\ell\uep(x,t)\vep(x,t),\qquad\qquad\qquad\qquad (x,t)\in \Om\times(0,T_{max,\e}),
\end{align*}
the Young inequality combined with \eqref{0vnav-inf}-\eqref{fu} yields that there exists $c_5(\e)>0$ such that
\begin{align*}
a_{\e}(x,t,\uep,\na\uep)\cdot\na\uep&=\vep\uep^{m-1}|\na\uep|^2
-f(\uep)\vep\na\vep\cdot\na\uep\no\\
&\ge c_3(\e)\uep^{m-1}|\na\uep|^2-c_4(\e)\uep^{\f{m-1}{2}}\vep|\na\vep||\na\uep|\no\\
&\ge c_3(\e)\uep^{m-1}|\na\uep|^2-c_2^2(\e)c_4(\e)\uep^{\f{m-1}{2}}|\na\uep|\no\\
&\ge \f{c_3}{2}(\e)\uep^{m-1}|\na\uep|^2- c_5(\e)
\end{align*}
as well as
\begin{align*}
&|a_{\e}(x,t,\uep,\na\uep)|\le c_2(\e)\uep^{m-1}|\na\uep|
+c_2^2(\e)c_4(\e)\uep^{\f{m-1}{2}},\qquad\qquad {\rm{and}}\no\\
&|b_{\e}(x,t)|\le \ell c_1(\e)c_2(\e)\qquad\qquad\qquad\qquad {\rm{for~all~}}(x,t)\in \Om\times(0,T_{max,\e}),
\end{align*}
which guarantee the {existence} of $\theta_1=\theta_1(\e), \theta_2=\theta_2(\e)\in(0,1)$ such that $\uep\in C^{\theta_1,\f{\theta_1}{2}}\left(\overline\Om\times[0,T_{max,\e}]\right)$ and
$\vep\in C^{2+\theta_2,1+\f{\theta_2}{2}}
\left(\overline\Om\times[\f{1}{4}T_{max,\e},T_{max,\e}]\right)$ provided by the H\"{o}lder estimates in \cite{holder1993} and parabolic Schauder theory in \cite{LSU1968}.

Now we arrange the first equation in \eqref{s1} to the following form
\beno
{\uep}_t= A_{\e}(x,t)\Del \uep+B_{\e}(x,t)\cdot\na\uep+D_{\e}(x,t)\uep\qquad x\in\Omega,\,t\in(0,T_{max,\e}),
\eeno
where
\begin{align*}
&A_{\e}(x,t):=u_{\e}^{m-1}(x,t)\vep(x,t),\\
&B_{\e}(x,t):=(m-1)u_{\e}^{m-2}(x,t)\vep(x,t)\na\uep(x,t)+u_{\e}^{m-1}(x,t)\na\vep(x,t)
-f'(\uep)\vep(x,t)\na\vep(x,t)\qquad{\rm{and}}\\
&D_{\e}(x,t):=-\f{f\left(\uep(x,t)\right)}{\uep(x,t)}\vep(x,t)\Del\vep(x,t)
-\f{f\left(\uep(x,t)\right)}{\uep(x,t)}|\na \vep(x,t)|^2+\ell\vep(x,t),\qquad(x,t)\in \Om\times(0,T_{max,\e}).
\end{align*}
By the assumptions on $m$, $f$ and $\al$, one can find some positive constant $c_6(\e)$ satisfying
\beno
\f{f\left(\uep(x,t)\right)}{\uep(x,t)}\le c_6(\e)\qmf(x,t)\in \Om\times(0,T_{max,\e}),
\eeno
which by means of the boundedness of $\vep$, $\na\vep$ and $\Del\vep$ entails that there exists $c_7(\e)>0$ such that
\beno
D_{\e}\ge -c_7(\e)  \qquad{\rm{in}}~~\Om\times(\f{1}{4}T_{max,\e},T_{max,\e}).
\eeno
Then by the comparison principle we get that
\beq\label{0u-up}
\uep\ge \inf_{x\in\Om} \uep\big(x,\f{1}{4}T_{max,\e}\big)e^{-c_7(\e)\cdot\f{3}{4}T_{max,\e}}  \qquad{\rm{in}}~~\Om\times(\f{1}{4}T_{max,\e},T_{max,\e}).
\eeq
Furthermore, applying the first order parabolic H\"{o}lder regularity theory (\cite{lieberman1987ampa}), there exists $\theta_3=\theta_3(\e)\in(0,1)$ such that
$\uep\in C^{1+\theta_3,\f{1+\theta_3}{2}}
\left(\overline\Om\times[\f{1}{2}T_{max,\e},T_{max,\e}]\right)$. In particular, we have
\beno
\|\uep(\cdot, t)\|_{W^{1,\infty}(\Omega)}\le c_8 \qmf t\in\big(\f{1}{2}T_{max,\e}, T_{max,\e}\big)
\eeno
with $c_8=c_8(\e)>0$. This together with \eqref{0vnav-inf}, \eqref{0v-up} and \eqref{0u-up} shows that \eqref{exten-1} fails, so that our assertion is identified. We thus obtain \eqref{exten}. Finally, \eqref{vin}-\eqref{uv1} can be derived easily by some basic integration computations and maximum principle.
\Ep

Frown now on, $\ell\ge0$ is fixed. Without further explicit mentioning, we assume $u_0$ and $v_0$ always fulfill \eqref{indata}, and accordingly let $(\uep,\vep)$ and $T_{max,\e}$ be as yielded by Lemma \ref{lelocal}.

Now we present two crucial inequalities from \cite[Lemma 3.4]{winkler2022dcdsb}, which will be utilized frequently in the following context.

\begin{Lemma}\label{lephi}
Let $q\ge2$. Then every $\varphi\in C^2(\overline\Om)$ fulfilling $\varphi>0$ in $\Om$ and $\f{\p\varphi}{\p\nu}=0$ on $\p\Om$ satisfies
\begin{align}\label{phi1}
\iom \f{|\na \varphi|^{q+2}}{\varphi^{q+1}}\le (q+\sqrt n)^2\iom \f{|\na \varphi|^{q-2}}{\varphi^{q-3}}|D^2\ln \varphi|^2
\end{align}
and
\begin{align}\label{phi2}
\iom \f{|\na \varphi|^{q-2}}{\varphi^{q-1}}|D^2\varphi|^2\le (q+\sqrt n+1)^2\iom \f{|\na \varphi|^{q-2}}{\varphi^{q-3}}|D^2\ln \varphi|^2.
\end{align}
\end{Lemma}

The next lemma is related to the type of functional $\iom v_{\e}^{-q+1}|\na \vep|^{q}$ for $q\ge2$, which was first proved in \cite[Lemma 3.3]{winkler2022dcdsb} for the general domain with boundary integral $\int_{\p\Om}v_{\e}^{-q+1}|\na\vep|^{q-2}\cdot\f{\p|\na\vep|^2}{\p\nu}$ arising. It was not feasible to treat this boundary integral term as in \cite[Lemma 3.5]{winkler2022dcdsb} because the term $\int_0^{\infty}\iom \vep(\cdot,t)$ is unexpected to appear later. So, as previously mentioned in \cite{li2022jde}, the convexity on $\Om$ in two-dimensional setting is imperative to allow us to disregard the boundary integral.
\begin{Lemma}\label{lenavq}
Let $\Om\subset\R^n$ ($n\le 2$) be a bounded convex domain and $q\ge2$. Then for all $t\in (0,T_{max,\e})$ and $\e\in(0,1)$, there exist $C_1(q)>0$ and $C_2>0$ such that
\begin{align}
&\f{d}{dt}\iom \f{|\na \vep|^{q}}{v_{\e}^{q-1}}+\f{q}{2}\iom \f{|\na \vep|^{q-2}}{v_{\e}^{q-3}}|D^2\ln\vep|^2+(q-1)^2\iom \uep\f{|\na \vep|^q}{v_{\e}^{q-1}}\le C_1(q)\iom \uep^{\f{q+2}{2}}v_{\e}\label{navq}
\end{align}
and
\begin{align}
&\f{d}{dt}\iom \f{|\na \vep|^{4}}{v_{\e}^{3}}+\iom \f{|\na \vep|^{2}}{v_{\e}}|D^2\ln\vep|^2+\iom \uep\f{|\na\vep|^4}{v_{\e}^3}\le C_2\iom \vep|\na\uep|^2. \label{nav4-v}
\end{align}
\end{Lemma}
\Bp
From Lemma 3.3 in \cite{winkler2022dcdsb}, we have
\begin{align}
&\f{d}{dt}\iom \f{|\na \vep|^{q}}{v_{\e}^{q-1}}+q\iom \f{|\na \vep|^{q-2}}{v_{\e}^{q-3}}|D^2\ln\vep|^2+(q-1)^2\iom \uep\f{|\na \vep|^q}{v_{\e}^{q-1}}\no\\
&\le q(q-2+\sqrt n)\iom \uep \f{|\na \vep|^{q-2}}{v_{\e}^{q-2}}|D^2\vep|+\f{q}{2}\int_{\p\Om} \f{|\na \vep|^{q-2}}{v_{\e}^{q-1}}\cdot \f{\p|\na\vep|^2}{\p\nu} \label{navq1}
\end{align}
for all $t\in (0,T_{max,\e})$ and $\e\in(0,1)$. Due to $\f{\p|\na\vep|^2}{\p\nu}\le 0$ on $\p\Om$ by convexity of $\Om$ (\cite{Lions1980}), we can {disregard} the last term in \eqref{navq1}. Using \eqref{phi1} and \eqref{phi2} and the Young inequality, the penultimate term in \eqref{navq1} could be estimated as follows
\begin{align*}
&q(q-2+\sqrt n)\iom \uep \f{|\na \vep|^{q-2}}{v_{\e}^{q-2}}|D^2\vep|\no\\
&\le\f{q}{4(q+\sqrt n+1)^2}\iom \f{|\na \vep|^{q-2}}{v_{\e}^{q-1}}|D^2\vep|^2+c_1\iom \uep^2\f{|\na \vep|^{q-2}}{v_{\e}^{q-3}}\no\\
&\le \f{q}{4}\iom \f{|\na \vep|^{q-2}}{v_{\e}^{q-3}}|D^2\ln \vep|^2+c_1\iom \uep^2\f{|\na \vep|^{q-2}}{v_{\e}^{q-3}}\no\\
&\le \f{q}{4}\iom \f{|\na \vep|^{q-2}}{v_{\e}^{q-3}}|D^2\ln \vep|^2+\f{q}{4(q+\sqrt n)^2}\iom \f{|\na \vep|^{q+2}}{v_{\e}^{q+1}}+c_2\iom \uep^{\f{q+2}{2}}v_{\e}\no\\
&\le \f{q}{2}\iom \f{|\na \vep|^{q-2}}{v_{\e}^{q-3}}|D^2\ln \vep|^2+c_2\iom \uep^{\f{q+2}{2}}v_{\e}  \qmf t\in (0,T_{max,\e})~{\rm{and}} ~\e\in(0,1)
\end{align*}
with positive constants $c_1=c_1(q)$ and $c_2=c_2(q)$, which combined \eqref{navq1} gives \eqref{navq}.

According to $(2.12)$ in \cite{li2022jde}, one has
\begin{align}
&\f{d}{dt}\iom \f{|\na \vep|^{4}}{v_{\e}^{3}}+4\iom \f{|\na \vep|^{2}}{v_{\e}}|D^2\ln\vep|^2+\iom\uep\f{|\na\vep|^4}{v_{\e}^3}\le -4\iom \f{|\na \vep|^{2}}{v_{\e}^{2}}(\na\uep\cdot\na\vep) \label{nav41}
\end{align}
for all $t\in (0,T_{max,\e})$ and $\e\in(0,1)$. Again using the Young inequality and \eqref{phi1} with $q=4$, implies the existence of $c_3>0$ such that
\begin{align*}
-4\iom \f{|\na \vep|^{2}}{v_{\e}^{2}}(\na\uep\cdot\na\vep)
&\le \f{3}{(4+\sqrt n)^2}\iom \f{|\na \vep|^{6}}{v_{\e}^{5}}+c_3\iom \vep|\na\uep|^2\no\\
&\le 3\iom \f{|\na \vep|^{2}}{v_{\e}}|D^2\ln\vep|^2+c_3\iom \vep|\na\uep|^2
\end{align*}
for all $t\in (0,T_{max,\e})$ and $\e\in(0,1)$. This inserting \eqref{nav41} ensures \eqref{nav4-v}.
\Ep

In the following, we show a boundedness result of a space-time integral derived by a basic calculation.
\begin{Lemma}\label{levnav}
Let $\Om\subset\R^n$ ($n\le 2$) be a bounded convex domain. For all $K>0$ there exists $C(K)>0$ with the property that whenever \eqref{indata1} or \eqref{indata2} holds, we have
\beq
\int_{0}^{T_{max,\e}}\iom  \vep|\na\vep|^2\le C(K)\label{vnav} \qmf \e\in(0,1).
\eeq
\end{Lemma}
\Bp
Testing the second equation in \eqref{s1} by $v_\e^2$ and integrating it on $\Om$ could yield this conclusion, where the boundedness of $\iom v_{0\e}^3$ ensured by \eqref{indata1} or \eqref{indata2} is necessary.\Ep

At the last of this section, we present the outcomes of some standard $L^p$ testing procedures applied to the first equation in \eqref{s1}.

\begin{Lemma}\label{leup}
For all $t\in (0,T_{max,\e})$ and $\e\in(0,1)$, the following estimations hold:
\begin{itemize}
\item[(i)] {If} $p>1$, then we have
\beq
\f{1}{p}\f{d}{dt}\iom  \uep^{p}+\f{p-1}{2}\iom \uep^{p+m-3}\vep|\na\uep|^2\le \f{p-1}{2}\iom \uep^{p-m-1}f^2(\uep)\vep|\na\vep|^2+\ell\iom \uep^p\vep;\label{upge1}
\eeq
\item[(ii)] if $0<p<1$, then we have
\beq
-\f{1}{p}\f{d}{dt}\iom  \uep^{p}+\f{1-p}{2}\iom \uep^{p+m-3}\vep|\na\uep|^2\le \f{1-p}{2}\iom \uep^{p-m-1}f^2(\uep)\vep|\na\vep|^2;\label{up01}
\eeq
\item[(iii)] if $p<0$, then we have
\beq
\f{d}{dt}\iom  \uep^{p}+\f{p(p-1)}{2}\iom \uep^{p+m-3}\vep|\na\uep|^2\le \f{p(p-1)}{2}\iom \uep^{p-m-1}f^2(\uep)\vep|\na\vep|^2;\label{uple0}
\eeq
\item[(iv)] if $m=2$, then we have
\beq
\f{d}{dt}\iom  \uep\ln \uep+\f{1}{2}\iom \vep|\na\uep|^2\le \f{1}{2}\iom \f{f^2(\uep)}{\uep^2}\vep|\na\vep|^2+\ell\iom \uep\vep+\ell\iom \uep\vep\ln \uep;\label{ulnu}
\eeq
\item[(v)] if $m=3$, then we have
\beq
-\f{d}{dt}\iom  \ln \uep+\f{1}{2}\iom \vep|\na\uep|^2\le \f{1}{2}\iom \f{f^2(\uep)}{\uep^4}\vep|\na\vep|^2.\label{lnu}
\eeq
\end{itemize}
\end{Lemma}
\Bp
If $p>1$ or $p<0$, we simply use the first equation in \eqref{s1} along with the boundary conditions, an integration by parts and the Young inequality to entail that
\begin{align}
\f{d}{dt}\iom \uep^{p}
&=p\iom \uep^{p-1}\nabla\cdot\left(u_{\e}^{m-1}{\vep}\nabla {\uep}\right)-p\iom \uep^{p-1}\nabla\cdot\big(f({\uep}){\vep}\nabla {\vep}\big)+p\ell\iom \uep^{p}\vep \no\\
&= -p(p-1)\iom \uep^{p+m-3}\vep|\na \uep|^2+p(p-1)\iom \uep^{p-2}f(\uep)\vep\na \uep\cdot\na\vep+p\ell\iom \uep^{p}\vep\no\\
&\le-\f{p(p-1)}{2} \iom \uep^{p+m-3}\vep|\na \uep|^2+\f{p(p-1)}{2}\iom \uep^{p-m-1}f^2(\uep)\vep|\na \vep|^2+p\ell\iom \uep^{p}\vep\label{up}
\end{align}
for all $t\in (0,T_{max,\e})$ and $\e\in(0,1)$. When $p>1$, dividing \eqref{up} by $p$ {establishes} \eqref{upge1}. When $p<0$, the negativity of the last term in \eqref{up} guarantees \eqref{uple0}.

Noting that $p(1-p)>0$ when $0<p<1$, similar to the procedures of getting \eqref{up}, we have
\begin{align*}
-\f{d}{dt}\iom \uep^{p}
&= -p(1-p)\iom \uep^{p+m-3}\vep|\na \uep|^2+p(1-p)\iom \uep^{p-2}f(\uep)\vep\na \uep\cdot\na\vep-p\ell\iom \uep^{p}\vep\no\\
&\le-\f{p(1-p)}{2} \iom \uep^{p+m-3}\vep|\na \uep|^2+\f{p(1-p)}{2}\iom \uep^{p-m-1}f^2(\uep)\vep|\na \vep|^2
\end{align*}
for all $t\in (0,T_{max,\e})$ and $\e\in(0,1)$, and then dividing this inequality by $p$, we get \eqref{up01}.

The proof of \eqref{ulnu} follows a similar approach to the proof of \eqref{upge1}, and the proof of \eqref{lnu} is analogous to that of \eqref{up01}. We omit the detailed proofs here.
\Ep

\section{Uniform boundedness of $\uep$ when $n=1$}
In this section, we first introduce a crucial functional inequality coming from \cite{li-winkler2022cpaa}, which is very useful to address terms of the form $\iom \uep^{\beta}\vep$ appearing in Lemma \ref{lecru} and Lemma \ref{leup5} in one-dimensional framework, and is inaccessible in two-dimensional case.
\begin{Lemma}\label{lefi1}
Let $\Om\subset\R$, $p\ge1$ and $r>1$. Then for all $\eta>0$ there exists $C(\eta,p,r)>0$ such that
\beq\label{fi1}
\big\|\varphi^{\f{p+1}{2}}\sqrt\psi\big\|_{L^{r}(\Om)}^2\leq \eta\iom\varphi^{p-1}\psi\varphi_x^2+\eta\cdot\bigg\{\iom\varphi^p\bigg\}
\cdot\iom\f{\varphi}{\psi} \psi_x^2+C(\eta,p,r)\cdot\bigg\{\iom\varphi\bigg\}^p
\cdot\iom\varphi\psi
\eeq
is valid for arbitrary nonnegative function $\varphi\in C^1(\overline\Om)$ and positive function $\psi\in C^1(\overline\Om)$.
\end{Lemma}

The next is about {an} a priori bound on the gradient of the second solution component in
\eqref{s1}. As a direct application of this, in light of Lemma \ref{lefi1}, the terms with the form $\iom \uep^{\gamma}\vep v_{\e x}^2$ appearing in Lemma \ref{leF} and Lemma \ref{leup5} could be dealt.
\begin{Lemma}\label{lenav4}
Let $\Om\subset\R$ be an open interval and $K>0$. Then there exists $C(K)>0$ with the property that whenever \eqref{indata1} holds, then
\beq\label{nav4}
\|{v}_{\e x}(\cdot,t)\|_{L^4(\Om)}\le C(K)\qquad for~all~ t\in (0,T_{max,\e})~{{and}} ~\e\in(0,1).
\eeq
\end{Lemma}
\Bp
From \eqref{u1}, we have
\beno
\|\uep(\cdot,t)\|_{L^1(\Om)}\le\iom u_{0\e}+\ell\iom v_{0\e}\le \|u_0\|_{L^1(\Om)} +|\Om|\big(1+\ell\cdot\|v_0\|_{L^{\infty}(\Om)}\big)
\eeno
for all $t\in (0,T_{max,\e})$ and $\e\in(0,1)$.
The variation-of-constants representation and well-known smoothing properties of Neumann heat semigroup $(e^{t\Delta})_{t\ge 0}$ (\cite[Lemma 1.3]{Winkler2010jde}) indicate \eqref{nav4}.\Ep

In the course of utilizing Lemma \ref{lefi1} to handle the integral $\iom \uep^2\vep$ arising from testing the second summand in \eqref{0F}, another unfavorable term $\iom \vep u_{\e x}^2$ will appear, which suggests us how to choose appropriate $p$ in \eqref{0F} at a first stage. The following is a more concrete outcome based on Lemma \ref{leup}.
\begin{Lemma}\label{leF}
Let $\Om\subset\R$ be an open interval. Then if one of the following cases holds:\\
(i) $1\le m<2$, $f$ fulfills \eqref{f} and \eqref{f1} with $m-1\le\al\le m$;\\
(ii) $2\le m<4$, $f$ fulfills \eqref{f} and \eqref{f2} with $m-1\le\al\le m$,\\
there exist positive constants $c$ and $C$ such that for each $\e\in(0,1)$, the function $\mathcal{F}_{\e }$ defined on $t\in (0,T_{max,\e})$ by letting
\beno
\mathcal{F}_{\e}(t):=\left\{
\begin{split}
&~~\iom \uep^{3-m}&& {{when}}\quad1\le m<2 \quad{{or}}\quad 3< m<4,\\
&~~\iom \uep\ln\uep  &&{{when}}\quad m=2,\\
&-\iom  \ln\uep      &&{{when}}\quad m=3,\\
&-\iom  \uep^{3-m}   &&{{when}}\quad2< m<3
\end{split}
\right.
\eeno
satisfies
\begin{align}
&\mathcal{F}_{\e }'(t)
+c\iom \vep u_{\e x}^2\le C\iom \uep^{2}\vep v_{\e x}^2+C\iom \uep^{2}\vep+C\iom \vep v_{\e x}^2+C\iom \uep\vep\label{F}
\end{align}
for all $t\in (0,T_{max,\e})$ and $\e\in(0,1)$.
\end{Lemma}
\Bp
We prove this assertion in five cases.\\
{\bf{Case 1: $1\le m<2$.}} Noting that $(1+s)^{2\al-2}\le c_1+c_1s^{2\al-2}$ is valid for all $s\ge0$ with some positive constant $c_1$, it yields by taking $p=3-m>1$ in \eqref{upge1} that
\begin{align}
&\f{d}{dt}\iom  \uep^{3-m}+\f{(2-m)(3-m)}{2}\iom \vep u_{\e x}^2\no\\
&\le \f{(2-m)(3-m)}{2}C_f^2\iom \uep^{4-2m}(\uep+1)^{2\al-2}\vep v_{\e x}^2+(3-m)\ell\iom \uep^{3-m}\vep\no\\
&\le c_1C_f^2\iom \uep^{4-2m}\vep v_{\e x}^2+c_1C_f^2\iom \uep^{2-2m+2\al}\vep v_{\e x}^2+2\ell\iom \uep^{3-m}\vep\label{F1-1}
\end{align}
for all $t\in (0,T_{max,\e})$ and $\e\in(0,1)$. Because $0<4-2m\le2$ due to $1\le m<2$ and $0\le2-2m+2\al\le2$ due to $m-1\le \al\le m$, using the Young inequality, the first two terms on the right side of \eqref{F1-1} could be estimated by
\begin{align*}
c_1C_f^2\iom \uep^{4-2m}\vep v_{\e x}^2+c_1C_f^2\iom \uep^{2-2m+2\al}\vep v_{\e x}^2\le 2c_1C_f^2\iom \uep^{2}\vep v_{\e x}^2+2c_1C_f^2\iom \vep v_{\e x}^2
\end{align*}
for all $t\in (0,T_{max,\e})$ and $\e\in(0,1)$. Another application of the Young inequality together with the fact $1<3-m\le 2$ shows that the last term in \eqref{F1-1} has the following estimation
\begin{align*}
2\ell\iom \uep^{3-m}\vep\le 2\ell\iom \uep^{2}\vep+2\ell\iom \uep\vep \qmf t\in (0,T_{max,\e})~{and} ~\e\in(0,1).
\end{align*}
Inserting the above two inequalities into \eqref{F1-1} shows that if we take $c:=\f{(2-m)(3-m)}{2}$ and  $C:=\max\{2c_1C_f^2,2\ell\}$, then \eqref{F} holds.\\
{\bf{Case 2: $ m=2$.}} Recalling \eqref{ulnu}, it follows from  the fact $\ln \uep\le \uep$ and the Young inequality in light of $0\le2\al-2\le2$ that for all $t\in (0,T_{max,\e})$ and $\e\in(0,1)$,
\begin{align*}
\f{d}{dt}\iom  \uep\ln \uep+\f{1}{2}\iom \vep v_{\e x}^2
&\le \f{1}{2}C_f^2\iom \uep^{2\al-2}\vep v_{\e x}^2+\ell\iom \uep\vep+\ell\iom \uep\vep\ln \uep\\
&\le \f{1}{2}C_f^2\iom \uep^{2}\vep v_{\e x}^2+\f{1}{2}C_f^2\iom \vep v_{\e x}^2+\ell\iom \uep\vep+\ell\iom u_{\e}^2\vep,
\end{align*}
which indicates \eqref{F} by taking $c:=\f{1}{2}$ and  $C:=\max\{\f{1}{2}C_f^2,\ell\}$.\\
{\bf{Case 3: $ 2<m<3$.}} Taking $p:=3-m\in (0,1)$ in \eqref{up01} and making use of the Young inequality once more with $0\le2-2m+2\al\le2$ due to $m-1\le \al\le m$, we get
\begin{align*}
-\f{d}{dt}\iom  \uep^{3-m}+\f{(m-2)(3-m)}{2}\iom \vep u_{\e x}^2
&\le \f{(m-2)(3-m)}{2}C_f^2\iom \uep^{2-2m+2\al}\vep v_{\e x}^2\\
&\le C_f^2\iom \uep^{2}\vep v_{\e x}^2+C_f^2\iom \vep v_{\e x}^2
\end{align*}
for all $t\in (0,T_{max,\e})$ and $\e\in(0,1)$. Then \eqref{Fnav2} is obtained by taking $c:=\f{(2-m)(3-m)}{2}$ and $C:=C_f^2$.\\
{\bf{Case 4: $ m=3$.}} It is easy to verify from \eqref{lnu} that if $m=3$, one has
\begin{align*}
-\f{d}{dt}\iom  \ln \uep+\f{1}{2}\iom \vep u_{\e x}^2
\le \f{1}{2}C_f^2\iom \uep^{2\al-4}\vep v_{\e x}^2\le \f{1}{2}C_f^2\iom \uep^{2}\vep v_{\e x}^2+\f{1}{2}C_f^2\iom \vep v_{\e x}^2
\end{align*}
for all $t\in (0,T_{max,\e})$ and $\e\in(0,1)$, where the last inequality is provided by the Young inequality with the help of $0\le 2\al-4\le2$ due to $2\le \al\le3$. This implies \eqref{Fnav2} by letting $c:=\f{1}{2}$ and $C:=\f{1}{2}C_f^2$.\\
{\bf{Case 5: $3< m<4$.}} The proof of this case can be derived in a very similar, even simpler, way as in case 1. We point out that the only difference is ``taking $p=3-m<0$ in \eqref{uple0}" instead of ``taking $p=3-m>1$ in \eqref{upge1}" when we use Lemma \ref{leup}. \Ep

Based on Lemmas \ref{lefi1}-\ref{leF}, we are able to derive the following conclusion by {making use of} the energy functional \eqref{0F}.
\begin{Lemma}\label{lecru}
Let $\Om\subset\R$ be an open interval, and $K>0$ with the property that \eqref{indata1} is valid. Then if one of the following cases holds:\\
(i) $1\le m<2$, $f$ fulfills \eqref{f} and \eqref{f1} with $m-1\le\al\le m$;\\
(ii) $2\le m<3$, $f$ fulfills \eqref{f} and \eqref{f2} with $m-1\le\al\le m$;\\
(iii) $3\le m<4$, $f$ fulfills \eqref{f} and \eqref{f2} with $m-1\le\al\le m$ and $u_0>0$ in $\overline{\Om}$,\\
one can find $C(K)>0$ such that
\beq\label{vux}
\int_0^{T_{max,\e}}\iom\uep \f{{v}_{\e x}^2}{\vep}+\int_0^{T_{max,\e}}\iom\vep u_{\e x}^2\le C(K) \quad\quad for~all~ \e\in(0,1).
\eeq
\end{Lemma}
\Bp
Taking $q=2$, \eqref{navq} in Lemma \ref{lenavq} has the following one-dimensional version
\beq\label{nav2}
\f{d}{dt}\iom \f{{v}_{\e x}^2}{\vep}+\iom \uep \f{{v}_{\e x}^2}{\vep}\le c_1\iom \uep^2v_{\e} \qmf t\in (0,T_{max,\e})~\rm{and}~\e\in(0,1)
\eeq
with $c_1>0$. Letting $\mathcal{F}_{\e}$ be as in Lemma \ref{leF}, it follows from \eqref{F} and \eqref{nav2} that there exist $c_2>0$ and $c_3>0$ fulfilling
\begin{align}
&\f{d}{dt}\iom\f{{v}_{\e x}^2}{\vep}+\mathcal{F}_{\e }'(t)
+c_2\iom \vep u_{\e x}^2+\iom\uep \f{{v}_{\e x}^2}{\vep}\no\\
&\le c_3\iom \uep^{2}\vep v_{\e x}^2+c_3\iom \uep^{2}\vep+c_3\iom \vep v_{\e x}^2+c_3\iom \uep\vep\label{Fnav2}
\end{align}
for all $t\in (0,T_{max,\e})$ and $\e\in(0,1)$, where the Cauchy-Schwarz inequality and \eqref{nav4} imply the existence of $c_4=c_4(K)>0$ such that
\begin{align}
c_3\iom \uep^{2}\vep v_{\e x}^2+c_3\iom \uep^{2}\vep&\le c_3\|v_{\e x}\|_{L^4(\Om)}^2 \cdot\|\uep{\vep}^{\f{1}{2}}\|_{L^4(\Om)}^2 +c_3|\Om|^{\f{1}{2}}\|\uep{\vep}^{\f{1}{2}}\|_{L^4(\Om)}^2\no\\
&\le c_4\|\uep{\vep}^{\f{1}{2}}\|_{L^4(\Om)}^2 \qmf t\in(0,T_{max,\e})~{\rm{and}}~\e\in(0,1),\label{11}
\end{align}
and we can further invoke Lemma \ref{lefi1} with $p:=1$ and $r:=4$ and use \eqref{u1} to find $c_5=c_5(K):=(K+1+\ell K)|\Om|$, $c_6=c_6(K)>0$ satisfying
\begin{align}
 c_4\|\uep{\vep}^{\f{1}{2}}\|_{L^4(\Om)}^2&\le
\f{c_2}{2}\iom {\vep}{u}_{\e x}^2+\f{1}{2c_5}\cdot\left\{\iom u_{\e}\right\}
\cdot\iom\uep \f{{v}_{\e x}^2}{\vep}+c_6\cdot\left\{\iom{\uep}\right\}
\cdot\iom{\uep}{\vep}\no\\
&\le \f{c_2}{2}\iom {\vep}{u}_{\e x}^2+\f{1}{2}\iom\uep \f{{v}_{\e x}^2}{\vep}+c_5c_6\iom{\uep}{\vep} \label{22}
\end{align}
for all $t\in(0,T_{max,\e})$ and $\e\in(0,1)$. Substituting \eqref{11} and \eqref{22} into \eqref{Fnav2}, we get
\begin{align}
&\f{d}{dt}\iom\f{{v}_{\e x}^2}{\vep}+\mathcal{F}_{\e }'(t)
+\f{c_2}{2}\iom \vep u_{\e x}^2+\f{1}{2}\iom\uep \f{{v}_{\e x}^2}{\vep}\le c_3\iom \vep v_{\e x}^2+(c_3+c_5c_6)\iom \uep\vep\label{Fnav2-1}
\end{align}
for all $t\in(0,T_{max,\e})$ and $\e\in(0,1)$.

A combination of \eqref{indata1} with Lemma \ref{levnav} and \eqref{uv1} entails that there exists $c_7=c_7(K)>0$ satisfying
\beno
c_3\int_0^{T_{max,\e}}\iom \vep v_{\e x}^2+
(c_3+c_5c_6)\int_0^{T_{max,\e}}\iom \uep\vep\le c_7 \qmf \e\in(0,1),
\eeno
as a consequence of which, upon an integration of \eqref{Fnav2-1}, we claim from \eqref{indata1} that
\begin{align*}
&\f{c_2}{2}\int_0^{t}\iom \vep u_{\e x}^2+\f{1}{2}\int_0^{t}\iom\uep \f{{v}_{\e x}^2}{\vep}\no\\
&\le c_3\int_0^{t}\iom \vep v_{\e x}^2+
(c_3+c_5c_6)\int_0^{t}\iom \uep\vep+\iom\f{{v}_{0 x}^2}{v_0}+ {\mathcal{F}_{\e }(0)}-{\mathcal{F}_{\e }(t)}\no\\
&\le c_7+K^3\cdot|\Om|+{\mathcal{F}_{\e }(0)}-{\mathcal{F}_{\e }(t)} \qmf t\in(0,T_{max,\e})~{\rm{and}}~\e\in(0,1).
\end{align*}
Thus, to complete the proof, it sufficies to show that ${\mathcal{F}_{\e }(0)}-{\mathcal{F}_{\e }(t)}$ always has an upper bound. Indeed, when $1\le m<2$, because of \eqref{indata1} and the positivity of $3-m$, we see that
\beno
{\mathcal{F}_{\e }(0)}-{\mathcal{F}_{\e }(t)}=\iom u_{0\e }^{3-m}-\iom u_{\e }^{3-m}\le \iom u_{0\e}^{3-m}\le(K+1)^{3-m}\cdot|\Om|.
\eeno
When $m=2$, the pointwise estimate $-\uep\ln \uep\le\f{1}{e}$ and \eqref{indata1} infer that
\beno
{\mathcal{F}_{\e }(0)}-{\mathcal{F}_{\e }(t)}=\iom u_{0\e}\ln u_{0\e} -\iom \uep\ln \uep\le (K+1)\ln (K+1)\cdot|\Om|+\f{|\Om|}{e}.
\eeno
When $2< m<3$, it follows from \eqref{u1} and the Young inequality with the fact $0<3-m<1$ that
\beno
{\mathcal{F}_{\e }(0)}-{\mathcal{F}_{\e }(t)}=-\iom u_{0\e}^{3-m}+\iom u_{\e }^{3-m}\le \iom u_{\e}^{3-m}\le \iom u_{\e}+|\Om|\le c_5+|\Om|.
\eeno
When $m=3$, the definition of $u_{0\e}$ and the strict positivity of $u_{0}$ entail the existence of $c_8>0$ independent on $\e$ fulfilling $u_{0\e}>c_8$. Thus, when $m=3$, due to the fact $\ln\uep\le\uep$, we have
\beno
{\mathcal{F}_{\e }(0)}-{\mathcal{F}_{\e }(t)}=-\iom \ln u_{0\e} +\iom \ln \uep\le-\iom \ln u_{0\e} +\iom \uep\le -\ln c_8\cdot|\Om|+c_5
\eeno
and when $3< m<4$, since $3-m$ is negative, we have
\beno
{\mathcal{F}_{\e }(0)}-{\mathcal{F}_{\e }(t)}=\iom u_{0\e}^{3-m}-\iom u_{\e }^{3-m}\le \iom u_{0\e}^{3-m}\le c_8^{3-m}\cdot|\Om|.
\eeno
The proof is finished. \Ep

We are now {prepared} to derive the main result of this section.
\begin{Lemma}\label{leup5}
Let $\Om\subset\R$ be an open interval, $p>2$ and $K>0$ with the property that \eqref{indata1} is valid. Then there exists $C(p,K)>0$ such that if one of the following cases holds:\\
(i) $1\le m<2$, $f$ fulfills \eqref{f} and \eqref{f1} with $m-1\le\al\le m$;\\
(ii) $2\le m<3$, $f$ fulfills \eqref{f} and \eqref{f2} with $m-1\le\al\le \f{m}{2}+1$;\\
(iii) $3\le m<4$, $f$ fulfills \eqref{f} and \eqref{f2} with $m-1\le\al\le \f{m}{2}+1$ and $u_0>0$ in $\overline{\Om}$,\\
for any choice of $\e\in(0,1)$, one has
\beq\label{eq-1}
\iom{u}_{\e}^p(\cdot,t)\le C(p,K)\quad for~all~ t\in (0,T_{max,\e}).
\eeq
\end{Lemma}
\Bp According to Lemma \ref{lenav4}, \eqref{u1} and \eqref{indata1}, we can fix $c_1=c_1(K)>0$ such that
\beq\label{nav4-1}
\iom{v}_{\e x}^4+\iom \uep\le c_1 \qmf t\in (0,T_{max,\e})~{\rm{and}}~\e\in(0,1),
\eeq
and we note that it is sufficient to prove \eqref{eq-1} for $p\ge4$ because of the Young inequality.

Now we firstly prove it in the case $1\le m<2$. From \eqref{upge1}, there exists $c_2=c_2(p)>0$ such that
\begin{align}
&\f{1}{p}\f{d}{dt}\iom  \uep^{p}+\f{p-1}{2}\iom \uep^{p+m-3}\vep u_{\e x}^2\no\\
&\le \f{p-1}{2}C_f^2\iom \uep^{p-m+1}(\uep+1)^{2\al-2}\vep v_{\e x}^2+\ell\iom \uep^p\vep\no\\
&\le c_2\iom \uep^{p-m-1+2\al}\vep v_{\e x}^2+c_2\iom \uep^{p-m+1}\vep v_{\e x}^2+\ell\iom \uep^p\vep\label{up5-1}
\end{align}
for all $t\in (0,T_{max,\e})$ and $\e\in(0,1)$. Making use of the Cauchy-Schwarz inequality, Lemma \ref{lefi1} and \eqref{nav4-1}, we can find $c_3=c_3(p,K)>0$ satisfying
\begin{align*}
c_2\iom \uep^{p-m-1+2\al}\vep v_{\e x}^2
&\le c_2\cdot\left\{\iom v_{\e x}^4\right\}^{\f{1}{2}}\cdot\left\{\iom \uep^{2(p-m-1+2\al)}v_{\e}^2\right\}^{\f{1}{2}}\\
&\le c_2c_1^{\f{1}{2}}\|{\uep}^{\f{p-m-1+2\al}{2}}v_{\e}^{\f{1}{2}}\|_{L^4(\Om)}^2\\
&\le \f{p-1}{4}\iom \uep^{p-m-3+2\al}\vep u_{\e x}^2+\left\{\iom \uep^{p-m-2+2\al}\right\}\cdot\iom\uep\f{{v}_{\e x}^2}{\vep}\\
&\qquad+c_3\cdot\left\{\iom{\uep}\right\}^{p-m-2+2\al}
\cdot\iom{\uep}{\vep}
\end{align*}
for all $t\in (0,T_{max,\e})$ and $\e\in(0,1)$, where $0\le p-m-3+2\al\le p+m-3$ due to $m-1\le\al\le m$ and $p\ge4$ enables us to use the Young inequality to estimate
\begin{align*}
&\iom \uep^{p-m-3+2\al}\vep u_{\e x}^2\le \iom \uep^{p+m-3}\vep u_{\e x}^2+\iom \vep u_{\e x}^2 \qmf t\in(0,T_{max,\e})~{\rm{and}}~\e\in(0,1)
\end{align*}
and another application of the Young inequality together with the fact $1\le p-m-2+2\al<p$ because of $m-1\le\al\le m$ implies that
\begin{align*}
&\iom \uep^{p-m-2+2\al}\le \iom \uep^{p}+|\Om| \qmf t\in(0,T_{max,\e})~{\rm{and}}~\e\in(0,1),
\end{align*}
which infer that with $c_4=|\Om|$ and $c_5=c_3c_1^{p-m-2+2\al}$, for all $t\in (0,T_{max,\e})$ and $\e\in(0,1)$ we have
\begin{align}
c_2\iom \uep^{p-m-1+2\al}\vep v_{\e x}^2
&\le \f{p-1}{4}\iom \uep^{p+m-3}\vep u_{\e x}^2+\f{p-1}{4}\iom \vep u_{\e x}^2\no\\
&\qquad+\left\{\iom \uep^{p}\right\}\cdot\iom\uep\f{{v}_{\e x}^2}{\vep}+c_4\iom\uep\f{{v}_{\e x}^2}{\vep}+c_5\iom{\uep}{\vep}.
\label{up5-2}
\end{align}
Similarly, thanks to $0<p-m+1,p+m-2<p\le p+m-1$ in the case $1\le m<2$, we could invoke the Young inequality several times and Lemma \ref{lefi1} once more to obtain that
\begin{align}
&c_2\iom \uep^{p-m+1}\vep v_{\e x}^2+\ell\iom \uep^p\vep\no\\
&\le c_2\iom \uep^{p+m-1}\vep v_{\e x}^2+c_2\iom \vep v_{\e x}^2+\ell\iom \uep^{p+m-1}\vep+\ell\iom \uep\vep\no\\
&\le\left(c_2c_1^{\f{1}{2}}+\ell|\Om|^{\f{1}{2}}\right)\cdot\|{\uep}^{\f{p+m-1}{2}}{\vep}^{\f{1}{2}}
\|_{L^4(\Om)}^2+c_2\iom \vep v_{\e x}^2+\ell\iom \uep\vep\no\\
&\le \f{p-1}{4}\iom \uep^{p+m-3}\vep u_{\e x}^2+\left\{\iom \uep^{p+m-2}\right\}\cdot\iom\uep\f{{v}_{\e x}^2}{\vep}\no\\
&\qquad+c_6\cdot\left\{\iom{\uep}\right\}^{p+m-2}
\cdot\iom{\uep}{\vep}+c_2\iom \vep v_{\e x}^2+\ell\iom \uep\vep\no\\
&\le \f{p-1}{4}\iom \uep^{p+m-3}\vep u_{\e x}^2+\left\{\iom \uep^{p}+c_4\right\}\cdot\iom\uep\f{{v}_{\e x}^2}{\vep}+c_2\iom \vep v_{\e x}^2+c_7\iom \uep\vep
\label{up5-3}
\end{align}
for all $t\in (0,T_{max,\e})$ and $\e\in(0,1)$ with $c_6=c_6(p)>0$ and $c_7=c_6c_1^{p+m-2}+\ell$. This together with \eqref{up5-1} and \eqref{up5-2} shows that
\begin{align*}
\f{1}{p}\f{d}{dt}\iom  \uep^{p}&\le \f{p-1}{4}\iom \vep u_{\e x}^2+2\left\{\iom \uep^{p}\right\}\cdot\iom\uep\f{{v}_{\e x}^2}{\vep}+2c_4\iom\uep\f{{v}_{\e x}^2}{\vep}\no\\
&~~~~+c_2\iom \vep v_{\e x}^2+(c_5+c_7)\iom \uep\vep \qmf t\in(0,T_{max,\e})~ {\rm{and}}~ \e\in(0,1).
\end{align*}
If we write
\beno
y_{\e}(t):=\iom \uep^{p}(\cdot,t), ~~t\in[0,T_{max,\e})
\eeno
and
\beno
g_{\e}(t):=2p\iom \uep(\cdot,t)\f{{v}_{\e x}^2(\cdot,t)}{\vep(\cdot,t)},~~t\in(0,T_{max,\e})
\eeno
as well as
\begin{align*}
h_{\e}(t):&=\f{p(p-1)}{4}\iom \vep(\cdot,t) u_{\e x}^2(\cdot,t)+2pc_4\iom\uep(\cdot,t)\f{{v}_{\e x}^2(\cdot,t)}{\vep(\cdot,t)}\no\\
&\qquad+pc_2\iom \vep(\cdot,t) v_{\e x}^2(\cdot,t)+p(c_5+c_7)\iom{\uep(\cdot,t)}{\vep(\cdot,t)},~~t\in(0,T_{max,\e}),
\end{align*}
then we get that
\beno
y_{\e}'(t)\le g_{\e}(t)y_{\e}'(t)+h_{\e}(t) \qmf t\in(0,T_{max,\e})~ {\rm{and}}~ \e\in(0,1)
\eeno
with the property that there is $c_8=c_8(p,K)>0$ such that
\beno
\int_{0}^{t}g_{\e}(s)ds\le c_8 \quad{\rm{and}}\quad \int_{0}^{t}h_{\e}(s)ds\le c_8 \qmf t\in(0,T_{max,\e})~ {\rm{and}}~ \e\in(0,1)
\eeno
 provided by Lemma \ref{levnav}, Lemma \ref{lecru}, \eqref{uv1} and \eqref{indata1}. An ODE comparison argument in conjunction with the fact that $y_{\e}(0)=\iom (u_0+\e)^p\le(K+1)^p|\Om|$ will indicate \eqref{eq-1}.

In the case $2\le m<4$, invoking $f(\uep)\le C_f u_{\e}^{\al}$ into \eqref{upge1}, then \eqref{up5-1} becomes
\begin{align*}
&\f{1}{p}\f{d}{dt}\iom  \uep^{p}+\f{p-1}{2}\iom \uep^{p+m-3}\vep u_{\e x}^2\le c_2\iom \uep^{p-m-1+2\al}\vep v_{\e x}^2+\ell\iom \uep^p\vep
\end{align*}
for all $t\in(0,T_{max,\e})$ and $\e\in(0,1)$, where on the basis of $\f{m}{2}+1\le m$, the term $c_2\iom \uep^{p-m-1+2\al}\vep v_{\e x}^2$ can be estimated in the same way used in obtaining \eqref{up5-2}, and we point out that the rightmost term $\ell\iom \uep^p\vep$ will be handled in a slightly different manner from \eqref{up5-3}. Indeed, if $2\le m<4$, we have $0<p-1\le p+m-3$, which allows us to use Lemma \ref{lefi1} and the Young inequality to find $c_9=c_9(p)>0$ such that for all $t\in(0,T_{max,\e})$ and $\e\in(0,1)$,
\begin{align*}
\ell\iom \uep^p\vep&\le \ell\iom \uep^{p+1}\vep+\ell\iom \uep\vep\\
&\le\ell|\Om|^{\f{1}{2}}\cdot\|{\uep}^{\f{p+1}{2}}{\vep}^{\f{1}{2}}
\|_{L^4(\Om)}^2+\ell\iom \uep\vep\\
&\le \f{p-1}{4}\iom \uep^{p-1}\vep u_{\e x}^2+\left\{\iom \uep^{p}\right\}\cdot\iom\uep\f{{v}_{\e x}^2}{\vep}+c_9\cdot\left\{\iom{\uep}\right\}^{p}
\cdot\iom{\uep}{\vep}+\ell\iom \uep\vep\\
&\le \f{p-1}{4}\iom \uep^{p+m-3}\vep u_{\e x}^2+\f{p-1}{4}\iom \vep u_{\e x}^2+\left\{\iom \uep^{p}\right\}\cdot\iom\uep\f{{v}_{\e x}^2}{\vep}+(c_9c_1^p+\ell)\iom{\uep}{\vep}.
\end{align*}
Trivially applying the procedures after \eqref{up5-3} thereby proves \eqref{eq-1} in the case $2\le m<4$. \Ep

\section{Uniform boundedness of $\uep$ when $n=2$}

\subsection{Some crucial estimations}
As shown in Section 3, the differential inequality \eqref{fi1} plays a critical role in helping us utilize the degenerate diffusive action of the first equation in \eqref{s1}. However, this approach does not work in two-dimensional space because the embedding $W^{1,1}(\Om)\hookrightarrow L^2(\Om)$ is only continuous, not compact. {To make sure that the term of the form $\iom\varphi^{p-1}\psi|\na \varphi|^2$ involves an arbitrarily small coefficient}, we adopt the following functional inequality (\cite[Lemma 3.1]{li2022jde}).
\begin{Lemma}\label{lefi2}
Let $\Om\subset\R^2$, $p>0$ and $r\ge 2$. Then for all $\eta>0$ there exists $C(\eta,p,r)>0$ such that
\beq\label{fi2}
\left\|\varphi^{\f{p+1}{2}}\sqrt\psi\right\|_{L^{r}(\Om)}^2\leq \eta\iom\varphi^{p-1}\psi|\na \varphi|^2+\eta\iom\varphi^{p+1}\f{|\na \psi|^2}{\psi} +C(\eta,p,r)\cdot\left\{\iom\varphi\right\}^p
\cdot\iom\varphi\psi
\eeq
is valid for arbitrary nonnegative function $\varphi\in C^1(\overline\Om)$ and positive function $\psi\in C^1(\overline\Om)$.
\end{Lemma}
\Bp When $p\ge1$, for any given $\eta>0$, we take
\beno
\eta_1\equiv\eta_1(\eta,p)\le \min\bigg\{\f{1}{1+|\Om|^{\f{r-2}{r}}},\f{2\eta}{(p+1)^2}\bigg\}.
\eeno
Because $W^{1,2}(\Om)$ is compactly embedded into $L^r(\Om)$ and the latter space is embedded into $L^1(\Om)$ with the continuous sense, the Ehrling lemma implies that for $\eta_1$, there exists a positive constant $C_1=C_1(\eta_1,r)>0$ such that
\begin{align}
\left\|\varphi^{\f{p+1}{2}}\sqrt\psi\right\|_{L^{r}(\Om)}^2&\leq \f{\eta_1}{8}\left\|\na\left(\varphi^{\f{p+1}{2}}\sqrt\psi\right)\right\|_{W^{1,2}(\Om)}^2
+C_1\left\|\varphi^{\f{p+1}{2}}\sqrt\psi\right\|_{L^{1}(\Om)}^2\no\\
&\leq\f{\eta_1}{4}\left\|\na\left(\varphi^{\f{p+1}{2}}\sqrt\psi\right)\right\|_{L^{2}(\Om)}^2
+\f{\eta_1}{4}\left\|\left(\varphi^{\f{p+1}{2}}\sqrt\psi\right)\right\|_{L^{2}(\Om)}^2
+C_1\left\|\varphi^{\f{p+1}{2}}\sqrt\psi\right\|_{L^{1}(\Om)}^2\no\\
&\leq\f{\eta_1}{4}\left\|\na\left(\varphi^{\f{p+1}{2}}\sqrt\psi\right)\right\|_{L^{2}(\Om)}^2
+\f{1}{4}\left\|\left(\varphi^{\f{p+1}{2}}\sqrt\psi\right)\right\|_{L^{r}(\Om)}^2
+C_1\left\|\varphi^{\f{p+1}{2}}\sqrt\psi\right\|_{L^{1}(\Om)}^2\label{eq41},
\end{align}
where the last inequality is due to the H\"{o}lder inequality. Then we can use the Young inequality to interpolate the norm in $L^1(\Om)$ between $L^r(\Om)$ and $L^{\f{2}{p+1}}(\Om)$, that is, with $C_2=C_2(\eta_1,p,r)>0$ we have
\beno
C_1\left\|\varphi^{\f{p+1}{2}}\sqrt\psi\right\|_{L^{1}(\Om)}^2\leq \f{1}{2}\left\|\varphi^{\f{p+1}{2}}\sqrt\psi\right\|_{L^{r}(\Om)}^2
+C_2\left\|\varphi^{\f{p+1}{2}}\sqrt\psi\right\|_{L^{\f{2}{p+1}}(\Om)}^2,
\eeno
which combined with \eqref{eq41} indicates that for any $p\ge1$, we have
\beno
\left\|\varphi^{\f{p+1}{2}}\sqrt\psi\right\|_{L^{r}(\Om)}^2\leq \eta_1\left\|\na\left(\varphi^{\f{p+1}{2}}\sqrt\psi\right)\right\|_{L^{2}(\Om)}^2
+4C_2\left\|\varphi^{\f{p+1}{2}}\sqrt\psi\right\|_{L^{\f{2}{p+1}}(\Om)}^2
\eeno
and thereby proves \eqref{fi2} by some basic calculus.

In the case when $0<p<1$, using the Ehrling lemma on the basis of the fact that $W^{1,2}(\Om)$ is compactly embedded into $L^r(\Om)$, while the latter space trivially lies in $L^{\f{2}{p+1}}(\Om)$, could directly derive the conclusion. \Ep

Parallel to that in Section 3, the key of our analysis in this section is also to control {the} unfavorable terms $\iom \uep^{\beta}\vep$ and $\iom \uep^{\gamma}\vep|\na \vep|^2$. Fortunately, an observation is that in view of Lemma \ref{lefi2}, this objective can be achieved when $\beta$ and $\gamma$ fall within the appropriate intermediate range, respectively. Additionally, we could catch some rationales why we take $p*$ in \eqref{0G} arbitrarily close to $3-m$, but not equal to the critical value here.
\begin{Lemma}\label{leu-v}
Let $\Om\subset\R^2$ and $K>0$ with the property that \eqref{indata2} is valid. If $\kappa\in(-1,0)$, then for any $\beta\in[1,\kappa+3)$ and $\eta>0$, there exists $C(K,\beta,\eta)>0$ such that
\beq\label{u-v}
\iom \uep^{\beta}\vep\le \eta\iom \uep^{\kappa}\vep|\na \uep|^2+\eta\iom \f{|\na \vep|^{2}}{v_{\e}}|D^2\ln\vep|^2+C(K,\beta,\eta)\iom \uep\vep
\eeq
for all $t\in (0,T_{max,\e})$ and $\e\in(0,1)$.
\end{Lemma}
\Bp We begin with establishing the validity of \eqref{u-v} for $\beta\in[\kappa+2,\kappa+3)$. According to \eqref{u1} and \eqref{indata2}, we can find $c_1=c_1(K)>0$ such that
\beq\label{4u1}
\|\uep(\cdot,t)\|_{L^1(\Om)}\le c_1 \qmf t\in (0,T_{max,\e})~{\rm{and}}~\e\in(0,1).
\eeq
Taking
\beno
\vartheta:=\f{1}{3+\kappa-\beta} \quad {\rm{and}}\quad \vartheta^*:=\f{\vartheta}{\vartheta-1}=\f{1}{\beta-\kappa-2},
\eeno
it follows from $\beta\in[\kappa+2,\kappa+3)$ that
\beno
\vartheta\ge 1 \quad {\rm{and}}\quad  \vartheta^*>1.
\eeno
Thus, for given $\eta>0$, we may use the H{\"o}lder inequality, Lemma \ref{lefi2} and \eqref{4u1} to conclude that there exists $c_2=c_2(K,\beta,\eta)>0$ satisfying
\begin{align}
\iom \uep^{\beta}\vep&=\iom\left(\uep^{\f{\kappa+2}{2}}\vep^{\f{1}{2}}\right)^2 \uep^{\beta-\kappa-2}\no\\
&\le\left\|\uep^{\f{\kappa+2}{2}}\vep^{\f{1}{2}}
\right\|_{L^{2\vartheta}(\Om)}^2\cdot\left\|{\uep}^{\beta-\kappa-2}
\right\|_{L^{\vartheta^*}(\Om)}\no\\
&\le(c_1^{\f{1}{\vartheta^*}}+1)\left\|\uep^{\f{\kappa+2}{2}}\vep^{\f{1}{2}}
\right\|_{L^{2\vartheta}(\Om)}^2\no\\
&\le \f{\eta}{2}\iom \uep^{\kappa}\vep|\na \uep|^2+\iom \uep^{\kappa+2} \f{|\na\vep|^2}{\vep}+c_2\iom \uep\vep\label{u-v1}
\end{align}
for all $t\in (0,T_{max,\e})$ and $\e\in(0,1)$. In treating the second term on the right side of \eqref{u-v1}, we first apply Lemma \ref{lenavq} with $q
:=4$ to see the existence of $c_3>0$ such that
\begin{align*}
\iom\f{|\na\vep|^6}{\vep^5}\le c_3\iom \f{|\na \vep|^{2}}{v_{\e}}|D^2\ln\vep|^2 \qmf t\in (0,T_{max,\e})~{\rm{and}}~\e\in(0,1).
\end{align*}
Then let
\beno
\chi:=-\f{2}{\kappa} \quad {\rm{and}}\quad \chi^*:=\f{\chi}{\chi-1}=\f{2}{\kappa+2},
\eeno
where our assumption $\kappa\in(-1,0)$ ensures that
\beno
\chi> 1 \quad {\rm{and}}\quad  \chi^*>1,
\eeno
whence by means of the Young inequality and the H{\"o}lder inequality, Lemma \ref{lefi2} and \eqref{4u1} are applied again so as to deduce that with $c_i=c_i(K,\beta,\eta)>0,~i=4,5$, we have
\begin{align*}
\iom \uep^{\kappa+2} \f{|\na\vep|^2}{\vep}
&\le \f{\eta}{2c_3}\iom\f{|\na\vep|^6}{\vep^5}+c_4\iom \uep^{\f{3}{2}\kappa+3}\vep\\
&\le\f{\eta}{2}\iom \f{|\na \vep|^{2}}{v_{\e}}|D^2\ln\vep|^2+
c_4\iom\left(\uep^{\f{\kappa+2}{2}}\vep^{\f{1}{2}}\right)^2
\uep^{\f{1}{2}\kappa+1}\\
&\le\f{\eta}{2}\iom \f{|\na \vep|^{2}}{v_{\e}}|D^2\ln\vep|^2+c_4\left\|\uep^{\f{\kappa+2}{2}}\vep^{\f{1}{2}}
\right\|_{L^{2\chi}(\Om)}^2\cdot\left\|\uep^{\f{1}{2}\kappa+1}
\right\|_{L^{\chi^*}(\Om)}\\
&\le\f{\eta}{2}\iom \f{|\na \vep|^{2}}{v_{\e}}|D^2\ln\vep|^2
+c_1^{\f{1}{\chi^*}}c_4\left\|\uep^{\f{\kappa+2}{2}}\vep^{\f{1}{2}}
\right\|_{L^{2\chi}(\Om)}^2\\
&\le \f{\eta}{2}\iom \f{|\na \vep|^{2}}{v_{\e}}|D^2\ln\vep|^2+\f{\eta}{4}\iom \uep^{\kappa}\vep|\na \uep|^2+\f{1}{2}\iom \uep^{\kappa+2} \f{|\na\vep|^2}{\vep}+c_5\iom \uep\vep
\end{align*}
for all $t\in (0,T_{max,\e})$ and $\e\in(0,1)$, which indicates that
\begin{align}
\iom \uep^{\kappa+2} \f{|\na\vep|^2}{\vep}
\le \eta\iom \f{|\na \vep|^{2}}{v_{\e}}|D^2\ln\vep|^2+\f{\eta}{2}\iom \uep^{\kappa}\vep|\na \uep|^2+2c_5\iom \uep\vep\label{u-v2}
\end{align}
for all $t\in (0,T_{max,\e})$ and $\e\in(0,1)$. Inserting \eqref{u-v2} into \eqref{u-v1} will yield \eqref{u-v}.

Now we illustrate that \eqref{u-v} also holds for $\beta\in[1,\kappa+2)$. Using the Young inequality, we see that
\beno
\iom  \uep^{\beta}\vep\le \iom \uep^{\kappa+2}\vep+\iom\uep\vep \qmf t\in (0,T_{max,\e})~{\rm{and}}~\e\in(0,1),
\eeno
which in conjunction with the fact that \eqref{u-v} holds for $\kappa+2$ completes the proof.\Ep

The following conclusion is a byproduct of the above lemma.
\begin{Lemma}\label{leuvnav}
Let $\Om\subset\R^2$ and $K>0$ with the property that \eqref{indata2} is valid. If $\kappa\in(-1,0)$, then for any $\ga\in[ 0,\f{\kappa}{2}+2)$ and $\eta>0$, there exists $C=C(K,\ga,\eta)>0$ such that
\begin{align}\label{uvnav}
\iom \uep^{\ga}\vep|\na \vep|^2
&\le \eta\iom \uep^{\kappa}\vep|\na \uep|^2+\eta\iom \uep \f{|\na\vep|^4}{\vep^3}+\eta\iom \f{|\na \vep|^{2}}{v_{\e}}|D^2\ln\vep|^2\no\\
&\qquad+\iom\vep|\na \vep|^2+C(K,\ga,\eta)\iom \uep\vep
\end{align}
for all $t\in (0,T_{max,\e})$ and $\e\in(0,1)$.
\end{Lemma}
\Bp If $\ga\in[ 0,1]$, using the Young inequality twice and relying on \eqref{indata2}, we see that
\begin{align*}
\iom \uep^{\ga}\vep|\na \vep|^2&\le \iom \uep\vep|\na \uep|^2+\iom \vep|\na \uep|^2\no\\
&\le \eta\iom \uep \f{|\na\vep|^4}{\vep^3}+\f{1}{4\eta}\iom \uep\vep^5+\iom \vep|\na \uep|^2\no\\
&\le \eta\iom \uep \f{|\na\vep|^4}{\vep^3}+\f{K^4}{4\eta}\iom \uep\vep+\iom \vep|\na \uep|^2
\end{align*}
for all $t\in (0,T_{max,\e})$ and $\e\in(0,1)$. This shows that \eqref{uvnav} holds for $\ga\in[ 0,1]$. On the other hand, when $\ga\in(1,\f{\kappa}{2}+2)$, it evident that $2\ga-1\in(1,\kappa+3)$. Therefore, it is available to make use of the Young inequality, \eqref{indata2} and Lemma \ref{leu-v} to see that with $c_1=c_1(K,\gamma,\eta)>0$, we have
\begin{align*}
\iom \uep^{\ga}\vep|\na \vep|^2
&\le \eta\iom \uep \f{|\na\vep|^4}{\vep^3}+\f{K^4}{4\eta}\iom \uep^{2\ga-1}\vep\no\\
&\le\eta\iom \uep^{\kappa}\vep|\na \uep|^2+\eta\iom \f{|\na \vep|^{2}}{v_{\e}}|D^2\ln\vep|^2+\eta\iom \uep \f{|\na\vep|^4}{\vep^3}+c_1\iom \uep\vep
\end{align*}
for all $t\in (0,T_{max,\e})$ and $\e\in(0,1)$. Thereby, we finish the proof.\Ep

\subsection{Uniform $L^{p_0}$ bounds on $u_{\e}$ for some $p_0>1$}

Let us first derive some basic information on the time evolution of the first two summands in \eqref{0G}.\begin{Lemma}\label{leG}
Let $\Om\subset\R^2$ be a bounded convex domain, and $K>0$ with the property that \eqref{indata2} is valid. Then if one of the following cases holds:\\
(i) $1\le m<2$, $f$ fulfills \eqref{f} and \eqref{f1} with $m-1<\al<m$;\\
(ii) $2\le m<4$, $f$ fulfills \eqref{f} and \eqref{f2} with $m-1<\al< m$,\\
one can find positive constants $c$ and $C$ such that for each $\e\in(0,1)$, the function $\mathcal{G}_{\e }$ defined on $t\in (0,T_{max,\e})$ by letting

\beno
\mathcal{G}_{\e}(t):=\left\{
\begin{split}
&~~c\iom  \uep^{3-m}+\iom\f{|\na\vep|^4}{\vep^3}&& {{when}}\quad1\le m<2 \quad{{or}}\quad 3< m<4,\\
&~~c\iom\uep\ln\uep+\iom\f{|\na\vep|^4}{\vep^3}  &&{{when}}\quad m=2,\\
&-c\iom  \ln\uep+\iom\f{|\na\vep|^4}{\vep^3}    &&{{when}}\quad m=3,\\
&-c\iom  \uep^{3-m}+\iom\f{|\na\vep|^4}{\vep^3}   &&{{when}}\quad2< m<3
\end{split}
\right.
\eeno
satisfies for all $t\in (0,T_{max,\e})$ and $\e\in(0,1)$,
\begin{align}
&\mathcal{G}_{\e }'(t)
+\iom \f{|\na \vep|^{2}}{v_{\e}}|D^2\ln\vep|^2+\iom\uep\f{|\na\vep|^4}{v_{\e}^3}+\iom \vep|\na\uep|^2\no\\
&\le  C\iom \uep^{2-2m+2\al}\vep|\na \vep|^2+C\iom \uep^{2}\vep+C\iom \uep\vep\quad when~~m=1\quad or ~~m=2,\label{G1}
\end{align}
\begin{align}
&\mathcal{G}_{\e }'(t)
+\iom \f{|\na \vep|^{2}}{v_{\e}}|D^2\ln\vep|^2+\iom\uep\f{|\na\vep|^4}{v_{\e}^3}+\iom \vep|\na\uep|^2\no\\
&\le  C\iom \uep^{2-2m+2\al}\vep|\na \vep|^2+C\iom \uep^{4-2m}\vep|\na \vep|^2+C\iom \uep^{3-m}\vep\quad when~~1<m<2\label{G2}
\end{align}
and
\begin{align}
&\mathcal{G}_{\e }'(t)+\iom \f{|\na \vep|^{2}}{v_{\e}}|D^2\ln\vep|^2+\iom\uep\f{|\na\vep|^4}{v_{\e}^3}+\iom \vep|\na\uep|^2\no\\
&\le  C\iom \uep^{2-2m+2\al}\vep|\na \vep|^2\quad when~~2<m<4\label{G3}.
\end{align}
\end{Lemma}
\Bp We first apply \eqref{nav4-v} in Lemma \ref{lenavq} to find $c_1>0$ such that
\begin{align}
&\f{d}{dt}\iom \f{|\na \vep|^{4}}{v_{\e}^{3}}+\iom \f{|\na \vep|^{2}}{v_{\e}}|D^2\ln\vep|^2+\iom\uep\f{|\na\vep|^4}{v_{\e}^3}\le c_1\iom \vep|\na\uep|^2 \label{vnav4}
\end{align}
for all $t\in (0,T_{max,\e})$ and $\e\in(0,1)$. When $1\le m<2$, taking $p:=3-m>1$ in \eqref{upge1} and dividing both sides by $\f{2-m}{2(c_1+1)}$, it follows that
\begin{align*}
&\f{2(c_1+1)}{(2-m)(3-m)}\f{d}{dt}\iom  \uep^{3-m}+(c_1+1)\iom \vep |\na \uep|^2\no\\
&\le (c_1+1)C_f^2\iom \uep^{4-2m}(\uep+1)^{2\al-2}\vep|\na \vep|^2+\f{2\ell(c_1+1)}{2-m}\iom \uep^{3-m}\vep
\end{align*}
for all $t\in (0,T_{max,\e})$ and $\e\in(0,1)$. Then if we write
\beno
\mathcal{G}_{\e}(t):=\f{2(c_1+1)}{(2-m)(3-m)}\iom  \uep^{3-m}+\iom\f{|\na\vep|^4}{\vep^3},
\eeno
the above inequality together with \eqref{vnav4} deduces that
\begin{align}
&\mathcal{G}_{\e }'(t)+\iom \f{|\na \vep|^{2}}{v_{\e}}|D^2\ln\vep|^2+\iom\uep\f{|\na\vep|^4}{v_{\e}^3}+\iom \vep |\na \uep|^2\no\\
&\le (c_1+1)C_f^2\iom \uep^{4-2m}(\uep+1)^{2\al-2}\vep|\na \vep|^2+\f{2\ell(c_1+1)}{2-m}\iom \uep^{3-m}\vep \label{G-1}
\end{align}
for all $t\in (0,T_{max,\e})$ and $\e\in(0,1)$. Obviously, \eqref{G2} can be derived directly by using the fact that $(s+1)^{2\al-2}\le c_2s^{2\al-2}+c_2$ holds for any $s\ge0$ with some $c_2>0$.

When $m=1$, the inequality $(s+1)^{2\al-2}\le s^{2\al-2}$ holds for any $s\ge 0$ due to our assumption $\al<1$, which in conjunction with \eqref{G-1} proves \eqref{G1} for the case $m=1$. Moreover, for this $\mathcal{G}_{\e }$, it is easy to verify from \eqref{vnav4} and \eqref{uple0} that \eqref{G3} holds for the case $3< m<4$.

When $m=2$, for the fixed $c_1$, we could get from \eqref{ulnu} that
\begin{align*}
&2(c_1+1)\f{d}{dt}\iom  \uep\ln \uep+(c_1+1)\iom \vep|\na\uep|^2\no\\
&\le (c_1+1)C_f^2\iom \uep^{2\al-2}\vep|\na\vep|^2+2(c_1+1)\ell\iom \uep\vep+2(c_1+1)\ell\iom \uep^2\vep
\end{align*}
for all $t\in (0,T_{max,\e})$ and $\e\in(0,1)$, which together with \eqref{vnav4} implies  \eqref{G1} for the case $m=2$ by taking
\beno
\mathcal{G}_{\e}(t):=2(c_1+1)\iom  \uep\ln\uep+\iom\f{|\na\vep|^4}{\vep^3}.
\eeno

Analogously, based on \eqref{vnav4}, we can prove \eqref{G3} for the case $2<m<3$ by letting
\beno
\mathcal{G}_{\e}(t):=-\f{2(c_1+1)}{(m-2)(3-m)}\iom  \uep^{3-m}+\iom\f{|\na\vep|^4}{\vep^3},
\eeno
taking $p:=3-m\in (0,1)$ in \eqref{up01} and dividing both sides by $\f{m-2}{2(c_1+1)}$. Relying on \eqref{lnu} and \eqref{vnav4}, we can present \eqref{G3} for the case $m=3$ with
\beno
\mathcal{G}_{\e}(t):=-2(c_1+1)\iom  \ln\uep+\iom\f{|\na\vep|^4}{\vep^3}.
\eeno
The proof is thereby finished. \Ep

From the above lemma, we are readily able to derive the $L^{p_0}$ bounds for $\uep$ with some $p_0>1$ in the case $1\le m<2$.
\begin{Lemma}\label{lecru2-1}
Let $\Om\subset\R^2$ be a bounded convex domain and $1\le m<2$. Suppose that $f$ fulfills \eqref{f} and \eqref{f1} with $m-1<\al<m$ and $K>0$ with the property that \eqref{indata2} is valid. Then there exist $p_0>1$ and $C(K,p_0)>0$ such that
\beno
\iom \uep^{p_0}(\cdot,t)\le C(K,p_0)\quad for~all~ t\in (0,T_{max,\e})~and ~\e\in(0,1)
\eeno
and
\beno
\int_0^{T_{\max,\e}}\iom\uep\f{|\na\vep|^4}{v_{\e}^3}+\int_0^{T_{\max,\e}}\iom \vep|\na\uep|^2\le C(K,p_0) \quad for~all~\e\in(0,1).
\eeno
\end{Lemma}
\Bp Let $\mathcal{G_{\e}}$ be defined as in Lemma \ref{leG}. We see that there is $c_1>0$ such that for all $t\in (0,T_{max,\e})$ and $\e\in(0,1)$, when $m=1$,
\begin{align}
&\mathcal{G}_{\e }'(t)
+\iom \f{|\na \vep|^{2}}{v_{\e}}|D^2\ln\vep|^2+\iom\uep\f{|\na\vep|^4}{v_{\e}^3}+\iom \vep|\na\uep|^2\no\\
&\le c_1\iom \uep^{2-2m+2\al}\vep|\na \vep|^2+c_1\iom \uep^{2}\vep+c_1\iom \uep\vep\label{G-1-1}
\end{align}
and when $1<m<2$,
\begin{align}
&\mathcal{G}_{\e }'(t)
+\iom \f{|\na \vep|^{2}}{v_{\e}}|D^2\ln\vep|^2+\iom\uep\f{|\na\vep|^4}{v_{\e}^3}+\iom \vep|\na\uep|^2\no\\
&\le c_1\iom \uep^{4-2m}\vep|\na \vep|^2+ c_1\iom \uep^{2-2m+2\al}\vep|\na \vep|^2+c_1\iom \uep^{3-m}\vep.\label{G-2}
\end{align}
Noting that $m-1<1<3-m$ due to $m<2$, and that $\max\{m+1-2\al, 3-5m+4\al\}<3-m$ due to our assumption $m-1<\al<m$, thus it is possible to pick $p_0>1$ satisfying
\beq\label{p0}
\max\{m+1-2\al, 3-5m+4\al\}<p_0<3-m,
\eeq
from which we can check that
\beq\label{p0-1}
0<p_0-m+2\al-1<2-2m+2\al<\f{p_0}{2}+\f{m}{2}+\f{1}{2}
\eeq
and moreover, together with the fact that $3-m\le3m-1$ because of $m\ge 1$, we have
\beq\label{p0-2}
0<p_0-m+1<\f{p_0}{2}+\f{m}{2}+\f{1}{2}.
\eeq
For the fixed $p_0>1$, we employ \eqref{upge1} to find $c_2>0$ fulfilling
\begin{align}
&\f{1}{p_0}\f{d}{dt}\iom  \uep^{p_0}+\f{p_0-1}{2}\iom \uep^{p_0+m-3}\vep|\na\uep|^2\no\\
&\le \f{p_0-1}{2}C_f^2\iom \uep^{p_0-m+1}(1+\uep)^{2\al-2}\vep|\na\vep|^2+\ell\iom \uep^{p_0}\vep\no\\
&\le c_2\iom \uep^{p_0-m+2\al-1}\vep|\na\vep|^2+c_2\iom \uep^{p_0-m+1}\vep|\na\vep|^2+\ell\iom \uep^{p_0}\vep\label{up-01}
\end{align}
for all $t\in (0,T_{max,\e})$ and $\e\in(0,1)$. Observing that our restrictions $1<p_0<3-m$ and $m\ge 1$ warrant $\kappa:=p_0+m-3\in (-1,0)$, thus we may invoke Lemma \ref{leuvnav} along with \eqref{p0-1} and \eqref{p0-2} to infer the existence of $c_3=c_3(K,p_0)>0$ such that
\begin{align}
&c_1\iom \uep^{2-2m+2\al}\vep|\na\vep|^2+c_2\iom \uep^{p_0-m+2\al-1}\vep|\na\vep|^2+c_2\iom \uep^{p_0-m+1}\vep|\na\vep|^2\no\\
&\le\f{p_0-1}{8}\iom \uep^{p_0+m-3}\vep|\na \uep|^2+\f{1}{2}\iom \uep \f{|\na\vep|^4}{\vep^3}+\f{1}{4}\iom\f{|\na \vep|^{2}}{v_{\e}}|D^2\ln\vep|^2\no\\
&\qquad+\iom\vep|\na \vep|^2+c_3\iom \uep\vep \qmf t\in (0,T_{max,\e})~{\rm{and}}~\e\in(0,1). \label{up02}
\end{align}
For the taken $\kappa$, it is clear that $1<p_0<\kappa+3$ and $1<3-m<\kappa+3$, which ensure that Lemma \ref{leu-v} is available to estimate
\begin{align}
&\ell\iom \uep^{p_0}\vep+c_1\iom \uep^{3-m}\vep\no\\
&\le\f{p_0-1}{8}\iom \uep^{p_0+m-3}\vep|\na \uep|^2+\f{1}{4}\iom\f{|\na \vep|^{2}}{v_{\e}}|D^2\ln\vep|^2+c_4\iom \uep\vep \label{up03}
\end{align}
for all $t\in (0,T_{max,\e})$ and $\e\in(0,1)$ with $c_4=c_4(K,p_0)>0$.

Combining \eqref{G-1-1}, \eqref{up-01}, \eqref{up02} and \eqref{up03}, it gives that when $m=1$,
\begin{align}
&\mathcal{G}_{\e }'(t)+\f{1}{p_0}\f{d}{dt}\iom  \uep^{p_0}+\f{1}{2}\iom\uep\f{|\na\vep|^4}{v_{\e}^3}+\iom \vep|\na\uep|^2\no\\
&\le \iom\vep|\na \vep|^2+(c_1+c_3+c_4)\iom \uep\vep\label{G-3}\qmf t\in (0,T_{max,\e})~{\rm{and}}~\e\in(0,1).
\end{align}
And we note that if $1<m<2$, the fact $7-5m<3-m$ warrants the possibility to pick $p_0$ not only satisfying \eqref{p0}, but also fulfilling
\beno
p_0>7-5m.
\eeno
Therefore, we readily find that
\beno
0<4-2m<\f{p_0}{2}+\f{m}{2}+\f{1}{2}.
\eeno
Another application of Lemma \ref{leuvnav} may furnish positive constant $c_5=c_5(K,p_0)>0$ such that
\begin{align*}
c_1\iom \uep^{4-2m}\vep|\na \vep|^2&\le\f{p_0-1}{8}\iom \uep^{p_0+m-3}\vep|\na \uep|^2+\f{1}{4}\iom\uep\f{|\na\vep|^4}{v_{\e}^3}+\f{1}{4}\iom\f{|\na \vep|^{2}}{v_{\e}}|D^2\ln\vep|^2\no\\
&\qquad+\iom\vep|\na \vep|^2+c_5\iom \uep\vep \qmf t\in (0,T_{max,\e})~{\rm{and}}~\e\in(0,1),
\end{align*}
which combined with \eqref{G-2}, \eqref{up-01}, \eqref{up02} and \eqref{up03} gives that when $1<m<2$,
\begin{align}
&\mathcal{G}_{\e }'(t)+\f{1}{p_0}\f{d}{dt}\iom  \uep^{p_0}+\f{1}{4}\iom\uep\f{|\na\vep|^4}{v_{\e}^3}+\iom \vep|\na\uep|^2\no\\
&\le 2\iom\vep|\na \vep|^2+(c_3+c_4+c_5)\iom \uep\vep \qmf t\in (0,T_{max,\e})~{\rm{and}}~\e\in(0,1).\label{G-4}
\end{align}
Integrating \eqref{G-3} and \eqref{G-4} in time respectively, Lemma \ref{levnav}, \eqref{uv1} and \eqref{indata2} show at once that whenever $1\le m<2$, with some $c_6=c_6(K,p_0)>0$ and $c_7=c_7(K,p_0)>0$, we have
\begin{align*}
&\f{1}{p_0}\iom \uep^{p_0}
+\f{1}{4}\int_0^t\iom\uep\f{|\na\vep|^4}{v_{\e}^3}+\int_0^t\iom \vep|\na\uep|^2\no\\
&\le 2\int_0^t\iom\vep|\na \vep|^2+(c_1+c_3+c_4+c_5)\int_0^t\iom \uep\vep+\mathcal{G}_{\e }(0)-\mathcal{G}_{\e }(t)\no\\
&\le c_{6}+c_7\iom (u_0+1)^{3-m}+\iom \f{|\na v_0|^{4}}{v_{0}^{3}}\no\\
&\le c_{6}+c_7(K+1)^{3-m}|\Om|+K^5|\Om|\qmf t\in (0,T_{max,\e})~{\rm{and}}~\e\in(0,1).
\end{align*}
Thereby the proof is finished.\Ep

Differing from the previous lemma, when $2\le m\le4$, $p^*$ belongs to $(-2,1)$ rather than $p^*>1$. Therefore, merely relying on functional \eqref{0G} is insufficient to obtain our desired estimate directly. For this, we additionally consider another component $\iom u^{p_0}$ with some $p_0>1$.
\begin{Lemma}\label{lecru2-2}
Let $\Om\subset\R^2$ be a bounded convex domain and $2\le m<3$. Suppose that $f$ fulfills \eqref{f} and \eqref{f2} with $m-1<\al<\f{m}{2}+1$ and $K>0$ with the property that \eqref{indata2} is valid. Then there exist $0<p_*<1$, $p_0>1$ and $C(K,p_*,p_0)>0$ such that
\beno
\iom \uep^{p_0}(\cdot,t)\le C(K,p_*,p_0)\quad for~all~t\in (0,T_{max,\e})~and ~\e\in(0,1)
\eeno
and
\beno
\int_0^{T_{\max,\e}}\iom\uep\f{|\na\vep|^4}{v_{\e}^3}+\int_0^{T_{\max,\e}}\iom \vep|\na\uep|^2\le C(K,p_*,p_0)\quad for~all~\e\in(0,1).
\eeno
\end{Lemma}
\Bp Firstly, we note that our assumption $m-1<\al<\f{m}{2}+1$ {makes it possible} to find $p_*>0$ satisfying
\beq\label{p*}
\max\{m+1-2\al, 4\al-3m-1\}<p_*<3-m\le 1.
\eeq
Evidently, we have $\f{p_*}{2}+\f{3}{2}m+\f{3}{2}-2\al>1$ and $ p_*+m>1$, so that it is also possible to fix $p_0$ such that
\beq\label{p0-21}
1<p_0<\min\bigg\{\f{p_*}{2}+\f{3}{2}m+\f{3}{2}-2\al, p_*+m\bigg\}.
\eeq
Taking $\kappa:=p_*+m-3$, it is easy to verify from \eqref{p*} that
 \beq\label{c1-0}
\kappa\in(-1,0).
\eeq
Moreover, in view of $\al>m-1$, it follows from \eqref{p0-21} that
\beq\label{c1-1}
0<p_0-m-1+2\al<\f{p_*}{2}+\f{m}{2}+\f{1}{2}=:\f{\kappa}{2}+2,
\eeq
and due to our choice of $p*$ in \eqref{p*}, we have
\beq\label{c1-2}
0<p_*-m-1+2\al<2-2m+2\al<\f{p_*}{2}+\f{m}{2}+\f{1}{2}=:\f{\kappa}{2}+2,
\eeq
where the last inequality uses that the restriction $m\ge2$ warrants that $p_*>4\al-3m-1\ge 3-5m+4\al$.

Next, let $\mathcal{G_{\e}}$ be defined as in Lemma \ref{leG}. From \eqref{G1} and \eqref{G3}, we can conclude that when $2\le m<3$, there is $c_1>0$ such that
\begin{align}
&\mathcal{G}_{\e }'(t)
+\iom \f{|\na \vep|^{2}}{v_{\e}}|D^2\ln\vep|^2+\iom\uep\f{|\na\vep|^4}{v_{\e}^3}+\iom \vep|\na\uep|^2\no\\
&\le c_1\iom \uep^{2-2m+2\al}\vep|\na \vep|^2+c_1\iom \uep^{2}\vep+c_1\iom \uep\vep\label{G12-1}
\end{align}
for all $t\in (0,T_{max,\e})$ and $\e\in(0,1)$. For the above fixed $p*$ and $p_0$, we draw on \eqref{up01} and \eqref{upge1} respectively to see that for all $t\in (0,T_{max,\e})$ and $\e\in(0,1)$, we have
\begin{align*}
-\f{1}{p_*}\f{d}{dt}\iom  \uep^{p_*}+\f{1-p_*}{2}\iom \uep^{p_*+m-3}\vep|\na\uep|^2
\le \f{1-p_*}{2}C_f^2\iom \uep^{p_*-m-1+2\al}\vep|\na\vep|^2
\end{align*}
and
\begin{align*}
&\f{1}{p_0}\f{d}{dt}\iom  \uep^{p_0}\le \f{p_0-1}{2}C_f^2\iom \uep^{p_0-m-1+2\al}\vep|\na\vep|^2+\ell\iom \uep^{p_0}\vep,
\end{align*}
which in conjunction with \eqref{G12-1} yield that
\begin{align}
&\mathcal{G}_{\e }'(t)
+\f{1}{p_0}\f{d}{dt}\iom  \uep^{p_0}-\f{1}{p_*}\f{d}{dt}\iom  \uep^{p_*}+\f{1-p_*}{2}\iom \uep^{p_*+m-3}\vep|\na\uep|^2\no\\
&~~~+\iom \f{|\na \vep|^{2}}{v_{\e}}|D^2\ln\vep|^2+\iom\uep\f{|\na\vep|^4}{v_{\e}^3}+\iom \vep|\na\uep|^2\no\\
&\le \f{1-p_*}{2}C_f^2\iom \uep^{p_*-m-1+2\al}\vep|\na\vep|^2+\f{p_0-1}{2}C_f^2\iom \uep^{p_0-m-1+2\al}\vep|\na\vep|^2\no\\
&~~~+c_1\iom \uep^{2-2m+2\al}\vep|\na \vep|^2+\ell\iom \uep^{p_0}\vep+c_1\iom \uep^{2}\vep+c_1\iom \uep\vep\label{G12-1-1}
\end{align}
for all $t\in (0,T_{max,\e})$ and $\e\in(0,1)$.
Now in view of \eqref{c1-0}, \eqref{c1-1} and \eqref{c1-2}, Lemma \ref{leuvnav} becomes applicable to detect that there is $c_2=c_2(K,p_*,p_0)>0$ fulfilling
\begin{align}
&\f{1-p_*}{2}C_f^2\iom \uep^{p_*-m-1+2\al}\vep|\na\vep|^2+\f{p_0-1}{2}C_f^2\iom \uep^{p_0-m-1+2\al}\vep|\na\vep|^2+c_1\iom \uep^{2-2m+2\al}\vep|\na \vep|^2\no\\
&\le\f{1-p_*}{4}\iom \uep^{p_*+m-3}\vep|\na \uep|^2+\f{1}{2}\iom \uep \f{|\na\vep|^4}{\vep^3}+\f{1}{2}\iom\f{|\na \vep|^{2}}{v_{\e}}|D^2\ln\vep|^2\no\\
&\qquad+\iom\vep|\na \vep|^2+c_2\iom \uep\vep \qmf t\in (0,T_{max,\e})~{\rm{and}}~\e\in(0,1). \label{G12-2}
\end{align}
On the other hand, it is obvious from \eqref{p0-21} that $1<p_0<p_*+m:=\kappa+3$ and $2<p_*+m:=\kappa+3$, which together with \eqref{c1-0} enable us to employ Lemma \ref{leu-v} to find $c_3=c_3(K,p_*,p_0)>0$ such that
\begin{align}
\ell\iom \uep^{p_0}\vep+c_1\iom \uep^{2}\vep\le\f{1-p_*}{4}\iom \uep^{p_*+m-3}\vep|\na \uep|^2+\f{1}{2}\iom\f{|\na \vep|^{2}}{v_{\e}}|D^2\ln\vep|^2+c_3\iom \uep\vep \label{G12-3}
\end{align}
for all $t\in (0,T_{max,\e})$ and $\e\in(0,1)$. Inserting \eqref{G12-2} and \eqref{G12-3} into \eqref{G12-1-1} {shows} that
\begin{align}
&\mathcal{G}_{\e }'(t)
+\f{1}{p_0}\f{d}{dt}\iom  \uep^{p_0}-\f{1}{p_*}\f{d}{dt}\iom \uep^{p_*}+\f{1}{2}\iom\uep\f{|\na\vep|^4}{v_{\e}^3}+\iom \vep|\na\uep|^2\no\\
&\le \iom\vep|\na \vep|^2+(c_1+c_2+c_3)\iom \uep\vep\qmf t\in (0,T_{max,\e})~{\rm{and}}~\e\in(0,1).\label{G12-4}
\end{align}
Before proceeding, similar to the assertions in the proof of Lemma \ref{lecru}, we claim that there is $c_4=c_4(K)>0$ such that ${\mathcal{G}_{\e }(0)}-{\mathcal{G}_{\e }(t)}\le c_4$ for all $t\in (0,T_{max,\e})$ and $\e\in(0,1)$. Indeed when $m=2$, with some positive constant $c_5>0$ we have
\begin{align*}
{\mathcal{G}_{\e }(0)}-{\mathcal{G}_{\e }(t)}&=c_5\iom u_{0\e}\ln u_{0\e}-c_5\iom \uep\ln \uep+\iom \f{|\na v_0|^{4}}{v_{0}^{3}}\le c_5(K+1)^{2}|\Om|+\f{c_5}{e}|\Om|+K^5|\Om|,
\end{align*}
 and when $2< m<3$, using \eqref{u1}, \eqref{indata2} and the Young inequality along with $0<3-m<1$, there exists $c_6>0$ such that
\begin{align*}
{\mathcal{G}_{\e }(0)}-{\mathcal{G}_{\e }(t)}\le c_6\iom u_{\e }^{3-m}+\iom \f{|\na v_0|^{4}}{v_{0}^{3}}\no\le c_6\bigg(\iom  \uep+|\Om|\bigg)+K^5|\Om|\le \big(c_6(\ell K+K+2)+K^5\big)|\Om|.
\end{align*}
{So, integrating \eqref{G12-4} from $0$ to $t$, by means of \eqref{uv1} and Lemma \ref{levnav}, once more using \eqref{u1}, \eqref{indata2} and the Young inequality along with the fact that $0<p_*<1$, we have
\begin{align*}
&\f{1}{p_0}\iom  \uep^{p_0}+\f{1}{2}\int_0^t\iom\uep\f{|\na\vep|^4}{v_{\e}^3}+\int_0^t\iom \vep|\na\uep|^2\no\\
&\le \int_0^t\iom\vep|\na \vep|^2+c_7\int_0^t\iom \uep\vep+\mathcal{G}_{\e }(0)-\mathcal{G}_{\e }(t)+\f{1}{p_*}\iom  \uep^{p_*}+\f{1}{p_0}\iom  u_{0\e}^{p_0}\no\\
&\le c_{8}+c_4+\f{\ell K+K+2}{p_*}|\Om|+\f{K^{p_0} }{p_0}|\Om|\qmf t\in (0,T_{max,\e})~{\rm{and}}~\e\in(0,1),
\end{align*}
where $c_7=c_1+c_2+c_3$ and $c_8=c_8(K)>0$.} The proof is complete.\Ep
\begin{Lemma}\label{lecru2-3}
Let $\Om\subset\R^2$ be a bounded convex domain, $3\le m<4$ and $K>0$. Suppose that $f$ fulfills \eqref{f} and \eqref{f2} with $m-1<\al<\f{m}{2}+1$. Then there exist $p_*<0$, $p_0>1$, $C_1(K,p_*,p_0)>0$ and $C_2(K,p_*)>0$ such that if not only \eqref{indata2} holds but also $u_0>0$ in $\overline{\Om}$, one has
\beno
\iom \uep^{p_0}(\cdot,t)\le C_1(K,p_*,p_0)\quad for~all~t\in (0,T_{max,\e})~and ~\e\in(0,1)
\eeno
and
\beno
\int_0^{T_{\max,\e}}\iom\uep\f{|\na\vep|^4}{v_{\e}^3}+\int_0^{T_{\max,\e}}\iom \vep|\na\uep|^2\le C_2(K,p_*)\quad for~all~\e\in(0,1).
\eeno
\end{Lemma}
\Bp Let $\mathcal{G_{\e}}$ be defined as in Lemma \ref{leG} with $c_1>0$ when $m=3$, and with $c_2>0$ when $3<m<4$. We firstly show that there exists some positive constant $c_3=c_3(K)$ such that
\beq\label{Gbbd}
{\mathcal{G}_{\e }(0)}-{\mathcal{G}_{\e }(t)}\le c_3 \qmf t\in (0,T_{max,\e})~{\rm{and}}~\e\in(0,1).
\eeq
Indeed, the strict positivity of $u_0$ in $\overline\Om$ warrants the existence of $c_4>0$ independent on $\e$ fulfilling $ u_{0\e}\ge c_4$. Then in light of \eqref{u1} and \eqref{indata2}, when $m=3$,
\begin{align*}
{\mathcal{G}_{\e }(0)}-{\mathcal{G}_{\e }(t)}&=-c_1\iom \ln u_{0\e} +c_1\iom \ln \uep+\iom \f{|\na v_0|^{4}}{v_{0}^{3}}\no\\
&\le -c_1\iom \ln u_{0\e} +c_1\iom \uep+\iom \f{|\na v_0|^{4}}{v_{0}^{3}}\no\\
&\le -c_1\ln c_4\cdot|\Om|+c_1(\ell K+K+1)|\Om|+K^5|\Om|,
\end{align*}
and when $3< m<4$, since $3-m$ is negative and $u_{\e}$ is nonnegative, we have
\beno
{\mathcal{G}_{\e }(0)}-{\mathcal{G}_{\e }(t)}=c_2\iom u_{0\e}^{3-m}-c_2\iom u_{\e }^{3-m}+\iom \f{|\na v_0|^{4}}{v_{0}^{3}}\le c_2c_4^{3-m}|\Om|+K^5|\Om|.
\eeno
Moreover, relying on \eqref{G3}, we can find $c_5>0$ {satisfying}
\begin{align}
&\mathcal{G}_{\e }'(t)
+\iom \f{|\na \vep|^{2}}{v_{\e}}|D^2\ln\vep|^2+\iom\uep\f{|\na\vep|^4}{v_{\e}^3}+\iom \vep|\na\uep|^2\le c_5\iom \uep^{2+2\al-2m}\vep|\na \vep|^2\label{G3-1}
\end{align}
for all $t\in (0,T_{max,\e})$ and $\e\in(0,1)$.

Then similar to the discussions in the proof of Lemma \ref{lecru2-2}, we can pick $p_*\in(2-m, 3-m)$ and $p_0>1$ such that the following three conditions
\beq\label{c2-0}
0<p_*+2\al-m-1<2+2\al-2m<\f{p_*}{2}+\f{m}{2}+\f{1}{2},
\eeq
and
\beq\label{c2-2}
0<p_0+2\al-m-1<\f{p_*}{2}+\f{m}{2}+\f{1}{2}=\f{p_*+m-3}{2}+2
\eeq
as well as
\beq\label{c2-3}
1<p_0<p_*+m=(p_*+m-3)+3
\eeq
are fulfilled simultaneously. For the fixed $p_*$, recalling \eqref{uple0}, we have
\begin{align*}
\f{d}{dt}\iom  \uep^{p_*}+\f{p_*(p_*-1)}{2}\iom \uep^{p_*+m-3}\vep|\na\uep|^2
\le \f{p_*(p_*-1)}{2}C_f^2\iom \uep^{p_*+2\al-m-1}\vep|\na\vep|^2
\end{align*}
for all $t\in (0,T_{max,\e})$ and $\e\in(0,1)$. This in conjunction with \eqref{G3-1} yields that
\begin{align*}
&\mathcal{G}_{\e }'(t)+
\f{d}{dt}\iom  \uep^{p_*}+\f{p_*(p_*-1)}{2}\iom \uep^{p_*+m-3}\vep|\na\uep|^2\no\\
&~~~+\iom \f{|\na \vep|^{2}}{v_{\e}}|D^2\ln\vep|^2+\iom\uep\f{|\na\vep|^4}{v_{\e}^3}+\iom \vep|\na\uep|^2\no\\
&\le c_5\iom \uep^{2+2\al-2m}\vep|\na \vep|^2+\f{p_*(p_*-1)}{2}C_f^2\iom \uep^{p_*+2\al-m-1}\vep|\na\vep|^2
\end{align*}
for all $t\in (0,T_{max,\e})$ and $\e\in(0,1)$, whereas thanks to \eqref{c2-0}, {applying Lemma \ref{leuvnav} will result in} that there is $c_6=c_6(K,p_*)>0$ fulfilling
\begin{align*}
&c_5\iom \uep^{2+2\al-2m}\vep|\na \vep|^2+\f{p_*(p_*-1)}{2}C_f^2\iom \uep^{p_*+2\al-m-1}\vep|\na\vep|^2\no\\
&\le\f{p_*(p_*-1)}{4}\iom \uep^{p_*+m-3}\vep|\na \uep|^2+\f{1}{2}\iom \uep \f{|\na\vep|^4}{\vep^3}+\f{1}{2}\iom\f{|\na \vep|^{2}}{v_{\e}}|D^2\ln\vep|^2\no\\
&\qquad+\iom\vep|\na \vep|^2+c_6\iom \uep\vep \qmf t\in (0,T_{max,\e})~{\rm{and}}~\e\in(0,1),
\end{align*}
so that
\begin{align*}
&\mathcal{G}_{\e }'(t)
+\f{d}{dt}\iom  \uep^{p_*}+\f{p_*(p_*-1)}{4}\iom \uep^{p_*+m-3}\vep|\na \uep|^2\no\\
&~~~+\f{1}{2}\iom \f{|\na \vep|^{2}}{v_{\e}}|D^2\ln\vep|^2+\f{1}{2}\iom\uep\f{|\na\vep|^4}{v_{\e}^3}+\iom \vep|\na\uep|^2\no\\
&\le \iom\vep|\na \vep|^2+c_6\iom \uep\vep \qmf t\in (0,T_{max,\e})~{\rm{and}}~\e\in(0,1).
\end{align*}
Upon integration in time, in view of Lemma \ref{levnav}, \eqref{uv1}, \eqref{Gbbd} and the fact that $p_*$ is negative, this implies that there is $c_7=c_7(K)>0$ satisfying
\begin{align}\label{tsin}
&\f{p_*(p_*-1)}{4}\int_0^t\iom \uep^{p_*+m-3}\vep|\na \uep|^2+\f{1}{2}\int_0^t\iom \f{|\na\vep|^{2}}{v_{\e}}|D^2\ln\vep|^2\no\\
&~~~+\f{1}{2}\int_0^t\iom\uep\f{|\na\vep|^4}{v_{\e}^3}+\int_0^t\iom \vep|\na\uep|^2+\int_0^t\iom\vep|\na \vep|^2+\int_0^t\iom \uep\vep\no\\
&\le 2\int_0^t\iom\vep|\na \vep|^2+(c_6+1)\int_0^t\iom \uep\vep+\mathcal{G}_{\e }(0)-\mathcal{G}_{\e }(t)+\iom u_{0\e }^{p_*} \no\\
&\le c_7+c_3+c_4^{p_*}|\Om|\qmf t\in (0,T_{max,\e})~{\rm{and}}~\e\in(0,1).
\end{align}

For the above $p_0$, it follows from \eqref{upge1} that for all $t\in (0,T_{max,\e})$ and $\e\in(0,1)$,
\beno
\f{1}{p_0}\f{d}{dt}\iom  \uep^{p_0}\le \f{p_0-1}{2}C_f^2\iom \uep^{p_0+2\al-m-1}\vep|\na\vep|^2+\ell\iom \uep^{p_0}\vep,
\eeno
where in light of \eqref{c2-2} and Lemma \ref{leuvnav}, we get
\begin{align*}
&\f{p_0-1}{2}C_f^2\iom \uep^{p_0+2\al-m-1}\vep|\na\vep|^2\\
&\le\iom \uep^{p_*+m-3}\vep|\na \uep|^2+\iom \uep \f{|\na\vep|^4}{\vep^3}+\iom\f{|\na \vep|^{2}}{v_{\e}}|D^2\ln\vep|^2\no\\
&~~~+\iom\vep|\na \vep|^2+c_8\iom \uep\vep \qmf t\in (0,T_{max,\e})~{\rm{and}}~\e\in(0,1)
\end{align*}
with $c_8=c_8(K,p_0)>0$, and in view of \eqref{c2-3}, Lemma \ref{leu-v} says that with $c_9=c_9(K,p_0)>0$ we have
\begin{align*}
&\ell\iom \uep^{p_0}\vep\le\iom \uep^{p_*+m-3}\vep|\na \uep|^2+\iom\f{|\na \vep|^{2}}{v_{\e}}|D^2\ln\vep|^2+c_9\iom \uep\vep
\end{align*}
for all $t\in (0,T_{max,\e})$ and $\e\in(0,1)$. We therefore obtain that
\begin{align*}
\f{1}{p_0}\f{d}{dt}\iom  \uep^{p_0}&\le 2\iom \uep^{p_*+m-3}\vep|\na \uep|^2+\iom \uep \f{|\na\vep|^4}{\vep^3}+2\iom\f{|\na \vep|^{2}}{v_{\e}}|D^2\ln\vep|^2\no\\
&~~~+\iom\vep|\na \vep|^2+(c_8+c_9)\iom \uep\vep \qmf t\in (0,T_{max,\e})~{\rm{and}}~\e\in(0,1),
\end{align*}
while once integrated in time, this together with \eqref{tsin} completes the lemma.\Ep

\subsection{Uniform $L^{p}$ bounds on $u_{\e}$ for any $p>1$}
Now we are readily to get the bounds for $u_{\e}$ with respect to the norm in $L^{p}(\Om)$ for any $p>1$. Before doing this, we first state the following conclusion, which will allow us to control some ill-contributions arising in Lemma \ref{leup1}.
\begin{Lemma}\label{leu-v-1}
Let $\Om\subset\R^2$, $m\ge1$ and $p>1$. Suppose that $K>0$ with the property that \eqref{indata2} holds, and that for some $p_0>1$ there exists $c(K,p_0)>0$ satisfying
\begin{align}\label{equp0}
\iom{u}_{\e}^{p_0}(\cdot,t)\le c(K,p_0) \quad for~all~t\in (0,T_{max,\e})~and ~\e\in(0,1).
\end{align}
Then for any $\beta\in[p+m-1,p_0+p+m-1)$ and $\eta>0$, there exists $C(K,p_0,\beta,\eta)>0$ such that
\beq\label{u-v-1}
\iom \uep^{\beta}\vep\le \eta\iom \uep^{p+m-3}\vep|\na \uep|^2+\eta\iom \f{|\na \vep|^{q-2}}{v_{\e}^{q-3}}|D^2\ln\vep|^2+C(K,p_0,\beta,\eta)\iom \uep\vep
\eeq
for all $t\in (0,T_{max,\e})$ and $\e\in(0,1)$, where $q>\f{2(p+m-1)}{p_0}$.
\end{Lemma}
\Bp We first observe that
\beq\label{1}
1<\f{(p+m-1)(q+2)}{q}<p_0+p+m-1
\eeq
due to our assumptions $q>\f{2(p+m-1)}{p_0}$, $p>1$ and $m\ge1$, and note that for any $\beta>1$ satisfying $\beta\in\left[p+m-1,\f{(p+m-1)(q+2)}{q}\right)$, the Young inequality implies that $\iom \uep^{\beta}\vep\le\iom \uep\vep+\iom \uep^{\f{(p+m-1)(q+2)}{q}}\vep$ for all $t\in (0,T_{max,\e})$ and $\e\in(0,1)$. Therefore we claim that to prove this lemma, it is sufficient to concentrate our focus on
\beq\label{c3-0}
\beta\in\left[\f{(p+m-1)(q+2)}{q}, p_0+p+m-1\right),
\eeq
which entails that
\beno
\vartheta^*:=\f{p_0}{\beta-p-m+1}>1 \quad {\rm{and}}\quad  \vartheta:=\f{\vartheta^*}{\vartheta^*-1}>1.
\eeno
Then for given $\eta>0$, using the H{\"o}lder inequality, Lemma \ref{lefi2} and \eqref{equp0}, there is $c_1=c_1(K,p_0,\beta,\eta)>0$ such that
\begin{align}
\iom \uep^{\beta}\vep&=\iom\left(\uep^{\f{p+m-1}{2}}\vep^{\f{1}{2}}\right)^2 \uep^{\beta-(p+m-1)}\no\\
&\le\left\|\uep^{\f{p+m-1}{2}}\vep^{\f{1}{2}}
\right\|_{L^{2\vartheta}(\Om)}^2\cdot\left\|{\uep}^{\beta-(p+m-1)}
\right\|_{L^{\vartheta^*}(\Om)}\no\\
&\le \f{\eta}{2}\iom \uep^{p+m-3}\vep|\na \uep|^2+\iom \uep^{p+m-1} \f{|\na\vep|^2}{\vep}+c_1\iom \uep\vep\label{u-v-11}
\end{align}
for all $t\in (0,T_{max,\e})$ and $\e\in(0,1)$. Here Lemma \ref{lephi} and two applications of the Young inequality along with the fact that $1<\f{(p+m-1)(q+2)}{q}<\beta$ from \eqref{1} and \eqref{c3-0} show that there exist $c_2=c_2(\eta)>0$ and $c_3=c_3(\eta,\beta)>0$ such that
\begin{align*}
\iom \uep^{p+m-1} \f{|\na\vep|^2}{\vep}
&\le \f{\eta}{2(q+\sqrt 2)^2}\iom \f{|\na\vep|^{q+2}}{\vep^{q+1}}+c_2\iom \uep^{\f{(p+m-1)(q+2)}{q}} \vep\no\\
&\le \f{\eta}{2}\iom \f{|\na \vep|^{q-2}}{v_{\e}^{q-3}}|D^2\ln\vep|^2+\f{1}{2}\iom \uep^{\beta} \vep+c_3\iom \uep \vep
\end{align*}
for all $t\in (0,T_{max,\e})$ and $\e\in(0,1)$, which inserted into \eqref{u-v-11} yields that
\beno
\iom \uep^{\beta}\vep\le \eta\iom \uep^{p+m-3}\vep|\na \uep|^2+\eta\iom \f{|\na \vep|^{q-2}}{v_{\e}^{q-3}}|D^2\ln\vep|^2+2(c_1+c_3)\iom \uep\vep
\eeno
for all $t\in (0,T_{max,\e})$ and $\e\in(0,1)$.\Ep

With the improved integrability of $\uep$ (Lemmas \ref{lecru2-1}-\ref{lecru2-3}) at hand, the strategy of further deriving the announced boundedness of $\uep$ relies on
analysing the energy functional with the form of \eqref{0G0} for conveniently large $p>1$ and appropriately chosen $q>2$. By suitably interpolating, some unfavorable summands with the expression $\iom \uep^{\beta}\vep$ for some $\beta$ will rise on the right sides of \eqref{up11} and \eqref{navq-v}. To cope with these terms, we hope that the corresponding parameters meet the assumption in Lemma \ref{leu-v-1}, so $q$ therein should be in an appropriate intermediate range. In these processes the condition $p_0>1$ is of critical importance.

\begin{Lemma}\label{leup1}
Let $\Om\in\R^2$ be a bounded convex domain and $K>0$ with the property that \eqref{indata2} is valid. Then for all $p>2$ there exists $C(K,p)>0$ such that if one of the following cases holds:\\
(i) $1\le m<2$, $f$ fulfills \eqref{f} and \eqref{f1} with $m-1<\al< m$;\\
(ii) $2\le m<3$, $f$ fulfills \eqref{f} and \eqref{f2} with $m-1<\al<\f{m}{2}+1$;\\
(iii) $3\le m<4$, $f$ fulfills \eqref{f} and \eqref{f2} with $m-1<\al<\f{m}{2}+1$ and $u_0>0$ in $\overline{\Om}$,\\
one has,
\beq\label{eq1}
\iom{u}_{\e}^p(\cdot,t)\le C(K,p) \quad for~all~t\in (0,T_{max,\e})~and ~\e\in(0,1),
\eeq
\beq\label{eq2}
\int_0^{T_{max,\e}}\iom\uep^{p+m-3}\vep|\na\uep|^2\le C(K,p)\quad for~all~\e\in(0,1)
\eeq
and
\beq\label{eq3}
\int_0^{T_{max,\e}}\iom\uep^{p+m-1}\vep\le C(K,p)\quad for~all~\e\in(0,1).
\eeq
\end{Lemma}
\Bp We first recall Lemmas \ref{lecru2-1}-\ref{lecru2-3} to see that there exist some $p_0>1$ and $c_1=c_1(K,p_0)>0$ such that
\beq\label{equp0-1}
\iom{u}_{\e}^{p_0}(\cdot,t)\le c_1\quad for~all~t\in (0,T_{max,\e})~and ~\e\in(0,1),
\eeq
and note that Lemma \ref{levnav}, \eqref{uv1} and \eqref{indata2} entail the existence of $c_2=c_2(K)>0$ fulfilling
\beq\label{c4-0}
\int_0^{T_{\max\e}}\iom \vep|\na\uep|^2+\int_0^{T_{\max\e}}\iom \uep\vep\le c_2 \qmf \e\in(0,1),
\eeq
which in view of the Young inequality implies that it is sufficient to prove \eqref{eq1}-\eqref{eq3} for $p\ge3$. As $\f{2(p+m-1)}{p_0}<2(p_0+p+m-2)$, we are able to pick $q:=q(p)>2(p+m-2)$
in a way that
\beno
q\in\left(\f{2(p+m-1)}{p_0},~2(p_0+p+m-2)\right),
\eeno
which clearly warrants
\beq\label{2}
p+m-1<\f{(p+m-1)(q+2)}{q}<p_0+p+m-1
\eeq
and
\beq\label{3}
p+m-1<\f{q+2}{2}<p_0+p+m-1.
\eeq

Our assumptions on $\al$ and $m$ yield that $0<p-m+1\le p+m-1$ and $0<p+2\al-m-1<p+m-1$. Thus, we may once more invoke the Young inequality several times to deduce from \eqref{upge1} and \eqref{phi1} that there is $c_3=c_3(p)>0$ such that
\begin{align}
&\f{1}{p}\f{d}{dt}\iom  \uep^{p}+\f{p-1}{2}\iom \uep^{p+m-3}\vep|\na\uep|^2\no\\
&\le \f{p-1}{2}\iom \uep^{p-m-1}f^2(\uep)\vep|\na\vep|^2+\ell\iom \uep^p\vep\no\\
&\le \f{p-1}{2}C_f^2\iom \uep^{p+2\al-m-1}\vep|\na\vep|^2+\f{p-1}{2}C_f^2\iom \uep^{p-m+1}\vep|\na\vep|^2+\ell\iom \uep^p\vep\no\\
&\le (p-1)C_f^2\iom \uep^{p+m-1}\vep|\na\vep|^2+(p-1)C_f^2\iom \vep|\na\vep|^2+\ell\iom \uep^{p+m-1}\vep+\ell\iom \uep\vep\no\\
&\le \f{q}{2(q+\sqrt2)^2}\iom \f{|\na \vep|^{q+2}}{v_{\e}^{q+1}}+c_3\iom \uep^{\f{(p+m-1)(q+2)}{q}}\vep^{\f{3q+4}{q}}+(p-1)C_f^2\iom \vep|\na\vep|^2\no\\
&\qquad+\ell\iom \uep^{p+m-1}\vep+\ell\iom \uep\vep\no\\
&\le\f{q}{2}\iom \f{|\na \vep|^{q-2}}{v_{\e}^{q-3}}|D^2\ln\vep|^2+c_3K^{\f{2q+4}{q}}\iom \uep^{\f{(p+m-1)(q+2)}{q}}\vep+(p-1)C_f^2\iom \vep|\na\vep|^2\no\\
&\qquad+\ell\iom \uep^{p+m-1}\vep+\ell\iom \uep\vep \qmf t\in (0,T_{max,\e})~{\rm{and}}~\e\in(0,1).\label{up11}
\end{align}
This in conjunction with the following inequality identified from Lemma \ref{lenavq}
\begin{align}\label{navq-v}
\f{d}{dt}\iom \f{|\na \vep|^{q}}{v_{\e}^{q-1}}+q\iom \f{|\na \vep|^{q-2}}{v_{\e}^{q-3}}|D^2\ln\vep|^2\le c_4\iom \uep^{\f{q+2}{2}}v_{\e}
\end{align}
for all $t\in (0,T_{max,\e})$ and $\e\in(0,1)$ with $c_4=c_4(p)>0$ implies that

\begin{align}
&\f{d}{dt}\left\{\f{1}{p}\iom  \uep^{p}+\f{|\na \vep|^{q}}{v_{\e}^{q-1}}\right\}+\f{p-1}{2}\iom \uep^{p+m-3}\vep|\na\uep|^2+\f{q}{2}\iom \f{|\na \vep|^{q-2}}{v_{\e}^{q-3}}|D^2\ln\vep|^2\no\\
&\le c_3K^{\f{2q+4}{q}}\iom \uep^{\f{(p+m-1)(q+2)}{q}}\vep+c_4\iom \uep^{\f{q+2}{2}}v_{\e}+\ell\iom \uep^{p+m-1}\vep+(p-1)C_f^2\iom \vep|\na\vep|^2\no\\
&\qquad+\ell\iom \uep\vep\qmf t\in (0,T_{max,\e})~{\rm{and}}~\e\in(0,1).\label{upvnav}
\end{align}
Owing to \eqref{equp0-1}, \eqref{2} and \eqref{3}, Lemma \ref{leu-v-1} becomes applicable so as to ensure that there exists $c_5=c_5(K,p)>0$ such that
\begin{align}\label{navq2}
&c_3K^{\f{2q+4}{q}}\iom \uep^{\f{(p+m-1)(q+2)}{q}}\vep+c_4\iom \uep^{\f{q+2}{2}}v_{\e}+\ell\iom \uep^{p+m-1}\vep\no\\
&\le \f{p-1}{4}\iom \uep^{p+m-3}\vep|\na \uep|^2+\f{q}{4}\iom \f{|\na \vep|^{q-2}}{v_{\e}^{q-3}}|D^2\ln\vep|^2+c_5\iom \uep\vep
\end{align}
for all $t\in (0,T_{max,\e})$ and $\e\in(0,1)$. Inserting this into \eqref{upvnav}, we get that
\begin{align*}
&\f{d}{dt}\left\{\f{1}{p}\iom  \uep^{p}+\f{|\na \vep|^{q}}{v_{\e}^{q-1}}\right\}+\f{p-1}{4}\iom \uep^{p+m-3}\vep|\na \uep|^2+\f{q}{4}\iom \f{|\na \vep|^{q-2}}{v_{\e}^{q-3}}|D^2\ln\vep|^2\no\\
&\le (p-1)C_f^2\iom \vep|\na\vep|^2+(\ell+c_5)\iom \uep\vep \qmf t\in (0,T_{max,\e})~{\rm{and}}~\e\in(0,1),
\end{align*}
which together with \eqref{c4-0} and \eqref{indata2}, upon an integration in time, entails that with $c_6=c_6(K,p)>0$, we have
\begin{align}\label{navq3}
&\f{1}{p}\iom  \uep^{p}+\f{p-1}{4}\int_0^t\iom \uep^{p+m-3}\vep|\na \uep|^2+\f{q}{4}\int_0^t\iom \f{|\na \vep|^{q-2}}{v_{\e}^{q-3}}|D^2\ln\vep|^2\le c_6
\end{align}
for all $t\in (0,T_{max,\e})$ and $\e\in(0,1)$. Thus \eqref{eq1} and \eqref{eq2} are included. Going back to integrating \eqref{navq2} in time and making use of \eqref{c4-0} and \eqref{navq3} conclude \eqref{eq3}.\Ep

\section{Proof of Theorem \ref{th1} and Theorem \ref{th2}}
Now by making use of Lemma \ref{leup5} and Lemma \ref{leup1}, we could derive the $W^{1,\infty}(\Om)$ boundedness of $\vep$ from the standard well-known smoothing properties of Neumann heat semigroup once more.
\begin{Lemma}\label{lenav-inf}
Suppose that the assumptions in Theorem \ref{th1} or Theorem \ref{th2} are satisfied with $K>0$. Then there exists $C(K)>0$ fulfilling
\beq\label{nav-inf}
\|\vep(\cdot,t)\|_{W^{1,\infty}(\Om)}\le C(K)\quad for~all~t\in (0,T_{max,\e})~and ~\e\in(0,1).
\eeq
\end{Lemma}
\Bp We rely on Lemma \ref{leup5} and Lemma \ref{leup1} to see that whenever $n=1$ or $n=2$, for any $p>n$ there is $c_1=c_1(K,p)>0$ such that
\beno
\|\uep(\cdot,t)\|_{L^p(\Om)}\le c_1 \qmf t\in (0,T_{max,\e})~{\rm{and}}~\e\in(0,1),
\eeno
which guarantees that \eqref{nav-inf} can be deduced in a manner that closely resembles the reasoning in Lemma \ref{lenav4}.\Ep

Having Lemmas \ref{leup5}, \ref{leup1} and \ref{lenav-inf} at hand, we can conclude by using a method very similar to \cite[Lemma 4.1]{winkler2022na} that the approximate solutions obtained in Lemma \ref{lelocal} are global. That is, we have the following lemma.
\begin{Lemma}\label{leglobal}
Suppose that the assumptions in Theorem \ref{th1} or Theorem \ref{th2} are satisfied. Then $T_{max,\e}=\infty$ for all $\e\in(0,1)$.
\end{Lemma}

In the following, we will show that in one-dimensional setting, by means of the elliptic Harnack type inequality documented in \cite {winkler-preprint}, the $L^{\infty}$ boundedness of $\uep$ can be derived by adapting a similar strategy in \cite{li-winkler2022cpaa}.
\begin{Lemma}\label{le-u-infty}
Suppose that the assumptions in Theorem \ref{th1} are satisfied with $K>0$ such that \eqref{indata1} holds. Then there exists $C(K)>0$ with the property that for all $\e\in(0,1)$ we have
\beq\label{u-infty}
\|\uep(\cdot,t)\|_{L^{\infty}(\Om)}\le C(K)\quad for~all~t>0.
\eeq
\end{Lemma}
\Bp This result can be proved in much the same way as Theorems 5.3 and 5.4 in \cite{li-winkler2022cpaa}. So we just show the sketch. Firstly, thanks to the Harnack type inequality
\beq\label{u-infty}
\vep(x,t)\ge c_1\|\vep(\cdot,t)\|_{L^{\infty}(\Om)}\qmf x\in \Omega~~t>0~~{\rm{and}}~~\e\in(0,1)
\eeq
with $c_1=c_1(K)>0$, it is easy to verify from \eqref{u1} and \eqref{uv1} that
\beno
L_{\e}:=\int_0^{\infty}\|\vep(\cdot,t)\|_{L^{\infty}(\Om)}dt, \qquad\e\in(0,1)
\eeno
is well-defined. Setting
\beno
\tau:=\phi_{\e}(t):=\f{1}{L_{\e}}\int_0^t\|\vep(\cdot,s)\|_{{L^{\infty}}(\Om)}ds,\qquad t\ge0
\eeno
and
\beno
w_{\e}(x,\tau):=\uep(x,\phi_{\e}^{-1}(\tau)),\qquad x\in\overline\Om,~\tau\in[0,1),
\eeno
it follows from \eqref{s1} that for each $\e\in(0,1)$
\beno
\left\{
\begin{split}
&w_{\e\tau}=\big(a_{\e}(x,\tau)w_{\e}^{m-1}w_{\e x}\big)_x-\big(b_{\e}(x,\tau)f(w_{\e})\big)_x+\ell a_{\e}(x,\tau)w_{\e},&&x\in\Omega,\,\tau\in(0,1),\\
&w_{\e x}=0, &&x\in\p\Omega,\,\tau\in(0,1),\\
&w_{\e}(x,0)=u_{0\e}(x), &&x\in\Omega
\end{split}
\right.
\eeno
is valid with
\beno
a_{\e}(x,\tau):=L_{\e}\cdot\f{\vep(x,t)}{\|\vep(\cdot,t)\|_{{L^{\infty}}(\Om)}} ~~~~{\rm{and}}~~~~b_{\e}(x,\tau):=L_{\e}\cdot\f{\vep(x,t)v_{\e x}(x,t)}{\|\vep(\cdot,t)\|_{{L^{\infty}}(\Om)}}.
\eeno
Additionally, there exists $c_2=c_2(K)>0$ such that
\beno
\f{1}{c_2}\le a_{\e}(x,\tau)\le c_2 \qmf x\in \Omega,~\tau\in(0,1)~~{\rm{and}}~~\e\in(0,1),
\eeno
and for any $p>1$ there is $c_3=c_3(p,K)>0$ satisfying
\beno
\|w_{\e}(\cdot,\tau)\|_{{L^{p}}(\Om)}+\|b_{\e}(\cdot,\tau)\|_{{L^{p}}(\Om)} \le c_3 \qmf \tau\in(0,1)~~{\rm{and}}~~\e\in(0,1),
\eeno
which combined with the Moser-type result (\cite[Lemma A.1]{tao-winkler2012jde1}) show the existence of $c_4=c_4(K)>0$ such that
\beno
\|w_{\e}(\cdot,\tau)\|_{{L^{\infty}}(\Om)}\le c_4 \qmf \tau\in(0,1)~~{\rm{and}}~~\e\in(0,1).
\eeno
By rescaling back to $\uep$, we finish the proof.\Ep

Actually, up to this point if we recall how the regularized form behaves in \eqref{s1} and \eqref{u0e}, we can assert that the proofs of case $(iii)$ in Theorem \ref{th1} and case $(iii)$ in Theorem \ref{th2} have already been finished. Therefore, in the rest of this context we just concentrate on the case $1\le m<3$.

When $n=1$, similar to the arguments in \cite[Theorems 5.4-5.7]{li-winkler2022cpaa}, we can proceed to identify a suitable sequence whose limit solution $(u,v)$ is accurately a weak solution of \eqref{s}, which simultaneously enjoys the additional properties announced in \eqref{th1r}.
\begin{Lemma}\label{le-se1}
Let $n=1$ and suppose that the assumptions in Theorem \ref{th1} are satisfied. Then there exist $(\e_j)_{j\in\N}\in(0,1)$ and functions $u\ge0$ and $v>0$ a.e. in $\Om\times(0,\infty)$ satisfying
\begin{equation*}
\left\{
\begin{split}
&u\in C^0(\overline\Om\times[0,\infty))\cap L^{\infty}(\Om\times(0,\infty))~~~~and\\
&v\in C^{2,1}(\overline\Om\times(0,\infty))
\end{split}
\right.
\end{equation*}
such that $\e=\e_j\searrow 0$ we have
\begin{align*}
&\uep\rightarrow u \qquad ~in~ C^{0}_{loc}(\overline\Om\times[0,\infty))\quad ~and \\
&\vep\rightarrow v\qquad  ~in~ C^{2,1}_{loc}(\overline\Om\times(0,\infty)),
\end{align*}
and that $(u,v)$ forms a global weak solution of \eqref{s} in the sense of Definition \ref{de}.
\end{Lemma}

When $n=2$, we first state some additional spatiotemporal regularity properties.
\begin{Lemma}\label{lest-1}
Suppose that the assumptions in Theorem \ref{th2} are satisfied with $K>0$ such that \eqref{indata2} holds. Then for any $p>4$, there exists $C_1(p,K)>0$ such that
\beq\label{eq4}
\int_0^t  \left\|\na\left(\uep^{\f{p+m-1}{2}}(\cdot,s)\vep(\cdot,s)
\right)\right\|_{L^{2}(\Om)}^2 ds \le C_1(K,p) \quad for~all~t>0~and~\e\in (0,1)
\eeq
and for any $T>0$, there exist $C_2(K,T,p)>0$ and $C_3(K,T)>0$ fulfilling
\beq\label{eq5}
\int_0^T  \left\|\p_t\left(\uep^{\f{p+m-1}{2}}(\cdot,t)\vep(\cdot,t)
\right)\right\|_{(W^{3,2}(\Om))^*}dt  \le C_2(K,T,p) \quad for~all~ \e\in (0,1)
\eeq
as well as
\beq\label{eq6}
\iom\ln\f{\|v_0\|_{L^{\infty}(\Om)}}{\vep(\cdot,t)}\le C_3(K,T) \quad for~all~ t\in (0,T)~and~\e\in (0,1).
\eeq
\end{Lemma}
\Bp Relying on \eqref{indata2} and \eqref{vin}, the Young inequality shows that
\begin{align*}
&\int_0^t \left\|\na\left(\uep^{\f{p+m-1}{2}}(\cdot,s)\vep(\cdot,s)
\right)\right\|_{L^{2}(\Om)}^2ds \no\\
&\le \f{(p+m-1)^2}{2}\int_0^t\iom \uep^{p+m-3}\vep^2 |\na\uep|^2+2\int_0^t\iom \uep^{p+m-1}|\na\vep|^2\no\\
&\le \f{(p+m-1)^2K}{2}\int_0^t\iom \uep^{p+m-3}\vep|\na\uep|^2+K^2\int_0^t\iom \uep^{2p+2m-3}\vep+\int_0^t\iom\uep\f{|\na\vep|^4}{v_{\e}^3}
\end{align*}
for all $t>0$ and $\e\in(0,1)$, which by means of \eqref{eq2}, \eqref{eq3} and Lemmas \ref{lecru2-1}-\ref{lecru2-3} indicates \eqref{eq4}.

To proceed to prove \eqref{eq5}, we first fix $\psi\in W^{3,2}(\Om)$ fulfilling $\|\psi\|_{W^{3,2}(\Om)}\le 1$. Drawing on the Sobolev inequality, one can find $c_1>0$ such that
\begin{align*}
\|\psi\|_{L^{\infty}(\Om)}+\|\na\psi\|_{L^{\infty}(\Om)}\le c_1  \qmf t>0~{\rm{and}}~\e\in(0,1).
\end{align*}
By computing directly and integrating by parts, we see that for all $t>0$ and $\e\in (0,1)$,
\begin{align}
&\iom {\p}_t \big(\uep^{\f{p+m-1}{2}}\vep
\big)\cdot\psi \no\\
&=-\f{p+m-1}{2}\iom \na\big(\uep^{\f{p+m-3}{2}}\vep\psi\big)
\cdot\left\{\uep^{m-1}\vep\na\uep-f(\uep)\vep\na\vep\right\}\no\\
&~~~+\f{p+m-1}{2}\ell\iom \uep^{\f{p+m-1}{2}}\vep^2\psi
-\iom \na\big(\uep^{\f{p+m-1}{2}}\psi\big)\cdot\na\vep-\iom \uep^{\f{p+m+1}{2}}\vep\psi\no\\
&=-\f{(p+m-1)(p+m-3)}{4}\bigg\{\iom \uep^{\f{p+3m-7}{2}}\vep^2|\na\uep|^2\psi-\iom \uep^{\f{p+m-5}{2}}f(\uep)\vep^2(\na\uep\cdot\na\vep)\psi\bigg\}\no\\
&~~~-\f{p+m-1}{2}\bigg\{\iom \uep^{\f{p+3m-5}{2}}\vep(\na\uep\cdot\na\vep)\psi+\iom \uep^{\f{p+3m-5}{2}}\vep^2(\na \uep\cdot\na\psi)-\ell\iom \uep^{\f{p+m-1}{2}}\vep^2\psi\no\\
&~~~-\iom \uep^{\f{p+m-3}{2}}f(\uep)\vep|\na\vep|^2\psi-\iom \uep^{\f{p+m-3}{2}}f(\uep)\vep^2(\na\vep\cdot\na\psi)+\iom \uep^{\f{p+m-3}{2}}(\na\uep\cdot\na\vep)\psi\bigg\}\no\\
&~~~-\iom \uep^{\f{p+m-1}{2}}(\na\vep\cdot\na\psi)-\iom \uep^{\f{p+m+1}{2}}\vep\psi\no\\
&=:-\f{(p+m-1)(p+m-3)}{4}\sum_{i=1}^2I_i(\e)-\f{p+m-1}{2}
\sum_{i=3}^8I_i(\e)+\sum_{i=9}^{10}I_i(\e).\label{I10}
\end{align}
Here, we see from $p>4$ and $1\le m<3$ that $0<\f{p+3m-7}{2}<p+m-3$ and $0\le 2m-2<p+m-1$ to assert that by means of the Young inequality if we pick $c_2=c_2(K)>0$ satisfying $\|\vep\|_{W^{1,\infty}(\Om)}\le c_2$, then for all $t>0$ and $\e\in(0,1)$,
\begin{align*}
|I_1|\le c_1c_2\iom \uep^{\f{p+3m-7}{2}}\vep|\na\uep|^2\le c_1c_2\iom \uep^{p+m-3}\vep|\na\uep|^2+c_1c_2\iom \vep|\na\uep|^2
\end{align*}
and
\begin{align*}
|I_3|+|I_4|&\le c_1c_2\iom \uep^{\f{p+3m-5}{2}}\vep|\na\uep|\\
&\le c_1c_2\iom \uep^{p+m-3}\vep|\na\uep|^2+c_1c_2\iom \uep^{2m-2}\vep\\
&\le c_1c_2\iom \uep^{p+m-3}\vep|\na\uep|^2+c_1c_2\iom \uep^{p+m-1}\vep+c_1c_2^2|\Om|.
\end{align*}
When $0\le\al<1$, we have $f(\uep)\le \uep$, which infers that for all $t>0$ and $\e\in(0,1)$,
\begin{align*}
|I_2|&\le c_1c_2^2\iom \uep^{\f{p+m-3}{2}}\vep|\na\uep|\le c_1c_2^2\iom \uep^{p+m-3}\vep|\na\uep|^2+c_1c_2^3|\Om|.
\end{align*}
When $\al\ge 1$, there exists $c_3>0$ such that $f(\uep)\le c_3\uep+c_3\uep^{\al}$. Then two applications of the Young inequality along with the facts that $0<\f{p+m-5}{2}+\al<p+m-3$ and $0\le2\al-2<p+m-1$ reveal that for all $t>0$ and $\e\in(0,1)$,
\begin{align*}
|I_2|&\le c_1c_2^2\iom \uep^{\f{p+m-5}{2}}f(\uep)\vep|\na\uep|\no\\
&\le c_1c_2^2c_3\iom \uep^{\f{p+m-3}{2}}\vep|\na\uep|+c_1c_2^2c_3\iom \uep^{\f{p+m-5}{2}+\al}\vep|\na\uep|\no\\
&\le 2c_1c_2^2c_3\iom \uep^{p+m-3}\vep|\na\uep|^2+c_1c_2^3c_3|\Om|+c_1c_2^2c_3\iom \uep^{2\al-2}\vep\\
&\le 2c_1c_2^2c_3\iom \uep^{p+m-3}\vep|\na\uep|^2+c_1c_2^2c_3\iom \uep^{p+m-1}\vep+2c_1c_2^3c_3|\Om|.
\end{align*}
Similarly, since $0<\f{p+m-3}{2}+\al<p+m-1$ and $0<\f{p+m-1}{2}<\f{p+m+1}{2}<p+m-1$, for all $t>0$ and $\e\in(0,1)$ we have
\begin{align*}
|I_6|+|I_7|&\le 2c_1c_2^2c_3\iom \uep^{\f{p+m-1}{2}}\vep+ 2c_1c_2^2c_3\iom \uep^{\f{p+m-3}{2}+\al}\vep\le 4c_1c_2^2c_3\iom \uep^{p+m-1}\vep+4c_1c_2^3c_3|\Om|
\end{align*}
and
\begin{align*}
|I_5|+|I_{10}|&\le c_1c_2\iom \uep^{\f{p+m-1}{2}}\vep+c_1\iom \uep^{\f{p+m+1}{2}}\vep\le (c_1c_2+c_1)\iom \uep^{p+m-1}\vep+(c_1c_2^2+c_1c_2)|\Om|.
\end{align*}
Using the Young inequality several times, we can treat $I_8$ and $I_9$ as follows
\begin{align*}
|I_8|+|I_9|&\le c_1\iom \uep^{\f{p+m-3}{2}}|\na \uep||\na \uep|+c_1\iom \uep^{\f{p+m-1}{2}}|\na \uep|\\
&\le c_1\iom \uep^{p+m-4}\vep|\na\uep|^2+c_1\iom \uep^{p+m-2}\vep+2c_1\iom \uep\f{|\na\vep|^2}{\vep}\\
&\le c_1\iom \uep^{p+m-3}\vep|\na\uep|^2+c_1\iom \uep^{p+m-1}\vep+2c_1\iom \uep\f{|\na\vep|^4}{\vep^3}+3c_1\iom\uep\vep+c_1\iom \vep|\na\uep|^2.
\end{align*}
Substituting these estimations on $I_i,~ i=1,2,\cdot\cdot\cdot,10$ into \eqref{I10}, we obtain that for all $t>0$ and $\e\in (0,1)$, there is $c_4=c_4(K)>0$ such that
\begin{align}
\left\|\p_t\big(\uep^{\f{p+m-1}{2}}\vep
\big)\right\|_{(W^{3,2}(\Om))^*} &\le c_4 \bigg\{\iom \uep^{p+m-3}\vep|\na\uep|^2+\iom \uep^{p+m-1}\vep\no\\
&~~~~+\iom\uep \f{|\na\vep|^4}{\vep^3}+\iom \vep|\na\uep|^2+\iom\uep\vep+1\bigg\}.\label{I10-2}
\end{align}
Integrating \eqref{I10-2} from $0$ to $T$, we obtain \eqref{eq5} as a consequence of \eqref{uv1}, \eqref{indata2}, Lemmas \ref{lecru2-1}-\ref{lecru2-3} and Lemma \ref{leup1}. Finally, \eqref{eq6} can be proved by using the same way in \cite[Lemma 4.4]{winkler2022na}. \Ep

With all the preparations, we are able to construct a weak solution to \eqref{s} when $n=2$. The presentation takes the following form. One can refer to \cite[Lemma 5.1]{winkler2022na} for the details.
\begin{Lemma}\label{lel}
Let $n=2$ and $p>2$. Suppose that the assumptions in Theorem \ref{th2} are satisfied. Then there exist $(\e_j)_{j\in\N}\in(0,1)$ and functions $u\ge0$ and $v>0$ a.e. in $\Om\times(0,\infty)$ satisfying
\begin{equation*}
\left\{
\begin{split}
&u\in L^{\infty}((0,\infty);L^{p}(\Om))\\
&v\in L^{\infty}(\Om\times(0,\infty))~~~~and\\
&\na v\in L^{\infty}(\Om\times(0,\infty))
\end{split}
\right.
\end{equation*}
such that $\e=\e_j\searrow 0$ we have
\begin{align*}
&\uep\rightarrow u \qquad a.e. ~in~ \Om\times(0,\infty)~and~in~ L_{loc}^{p}(\overline\Om\times[0,\infty)),\\
&\vep\rightarrow v\qquad a.e. ~in~ \Om\times(0,\infty)~ and ~in~ L_{loc}^{p}(\overline\Om\times[0,\infty)) \quad ~and \\
&\na\vep\stackrel{*}\rightharpoonup \na v\qquad  in~ L^{\infty}(\Om\times(0,\infty)),
\end{align*}
and that $(u,v)$ forms a global weak solution of \eqref{s} in the sense of Definition \ref{de}.
\end{Lemma}

Our main results thereby in fact reduce to a mere summary:\\
{\bf Proofs of Theorems \ref{th1} and \ref{th2}}. The proofs of cases $(i)$ and $(ii)$ in Theorem \ref{th1} are consequences of Lemma \ref{le-se1}, while Lemmas \ref{lenav-inf}-\ref{le-u-infty} together with Lemma \ref{lelocal} prove case $(iii)$ in Theorem \ref{th1}. Apart from that Lemma \ref{lel} indicates cases $(i)$ and $(ii)$ in Theorem \ref{th2}, and Lemmas \ref{leup1},\ref{lenav-inf},\ref{leglobal} in conjunction with Lemma \ref{lelocal} complete the proof of case $(iii)$ in Theorem \ref{th2}. \qquad $\Box$

\textbf{Acknowledgements.} The author is very grateful to Prof. Michael Winkler and Dr. Genglin Li for their useful discussions and advice. This work was supported by the China Scholarship Council (No. 202206290145) and the Deutsche Forschungsgemeinschaft (No. 462888149).

\small{
}

\end{document}